%Paper: alg-geom/9507013
%From: henri@math.uic.edu
%Date: Thu, 20 Jul 1995 21:35:15 -0500

  \input amstex
\documentstyle{amsppt}
\magnification=1200

\NoBlackBoxes
\topmatter

\title Descent, Motives and $K$-theory \endtitle

\rightheadtext{Descent, Motives and $K$-theory }

\author H. Gillet and C. Soul\'e \endauthor

\address Department of Mathematics, Statistics, and Computer Science,
(m/c249), University of Illinois at Chicago, 851 S. Morgan Street,
Chicago, IL 60607, U.S.A.\endaddress
\email henri\@math.uic.edu \endemail

\address Institut des Hautes \'{E}tudes Scientifiques,  35, Route de
Chartres, 91440, Bures-sur-Yvettes, France\endaddress

\thanks The first author was partially supported by N.S.F grant
DMS-9203379 \endthanks \endtopmatter \par \bigskip

In this paper we show how to write any variety over a field of
characteristic zero as the difference of two pure motives, thus
answering a question asked by Serre (in arbitrary characteristic,
\cite
{Se} p.341).

\medskip

The Grothendieck theory of (pure effective Chow) motives starts with
the category ${\bold V}$ of smooth projective varieties over a field
$k$. A {\it motive} is a pair $(X,p)$, where $X$ is in ${\bold V}$ and
$p$ is a projector in the ring of algebraic correspondences from $X$
to itself (algebraic cycles on $X
\times X$ modulo linear equivalence). One can add two motives
 $M_1 = (X_1 ,p_1)$ and $M_2 = (X_2 ,p_2)$ by considering
 $M_1 \oplus M_2 = \left( X_1 \amalg
X_2 \ ,
\ p_1 + p_2 \right)$. Let $K_0 ({\bold M})$ be the abelian group
associated to this monoid. Assuming the field $k$ has characteristic
zero, for any variety $X$ over $k$ (i.e. any reduced scheme of finite
type over $k$) we define a class $[X]$ in $K_0 ({\bold M})$, which is
characterized by the following two properties: if $X$ is connected and
lies in ${\bold V}$, $[X]$ is the class of the motive $(X,\Delta_X)$,
where $\Delta_X$ is the diagonal in $X \times X$; if $Y$ is a closed
subset in $X$, with its reduced scheme structure, the following
identity holds in $K_0 ({\bold M})$: $$ [X] = [Y] + [X-Y] . \leqno
(0.1) $$

\medskip

This result (Theorem 4) is derived from a stronger one. Namely, to any
variety $X$ we associate a cochain complex $W(X)$ in the additive
category of complexes of motives, which is well-defined up to
canonical homotopy equivalence. Our construction of $W(X)$ and the
proof of its properties (Theorem 2) use higher algebraic $K$-theory
and, more specifically, the {\it Gersten complexes} of schemes. These
complexes are made out of the $K$-theory of all residue fields of a
given scheme, and among their homology groups are precisely the Chow
groups of algebraic cycles modulo linear equivalence (see Section
1.1). To prove Theorem 2, we first extend Manin's identity principle
\cite{M1}, by showing that a complex $C.$ of varieties is contractible
as a complex of motives if and only if a certain family of Gersten
complexes associated to $C.$ are acyclic (Theorem 1).  We then use a
variant of the theory of cohomological descent which applies to
Gersten complexes (and to $K'$-theory) of simplicial schemes, and was
developped by the first author in
\cite{G2} (Proposition 1).
 Here the proper surjective maps used in the conventional theory of
cohomological descent (\cite{SD},
\cite{D}) are replaced by {\it envelopes} in the sense of \cite{F-G}
and \cite{G2}: a proper map of schemes $f:X' \rightarrow X$ is an
envelope when, for every field $F$, all $F$-valued points of $X$ lift
to $X'$.  When $X$ is a variety over a field $k$ of characteristic
zero, Chow's lemma and Hironaka's theorem imply that $X$ admits such
an envelope $X'$ which is both smooth and projective (resolution of
singularities is the only reason for us to ask that the ground field
has characteristic zero). This fact, together with simplicial
arguments similar to those of Deligne's theory of mixed Hodge
structures \cite{D}, lead to the existence of the complex $W(X)$. This
complex $W(X)$ is homotopy equivalent to a bounded complex, and its
class in $K_0 ({\bold M})$ is the virtual motive $[X]$. The proof of
(0.1) follows (using Theorem 1 and descent) from the fact that the
Gersten complexes of $Y$, $X$ and $X-Y$ fit into short exact
sequences, since points of $X$ lie either in $Y$ or in $X-Y$.

\medskip

We want to emphasize that we only consider pure motives, and that
mixed motives never enter our discussion. This is made possible by
Theorem 1, which enables us to avoid the usual difficulties arising
from the homotopy theory of complexes by giving us an acyclicity
criterion for a complex of motives to be contractible. If $k={\Bbb C}$
say, one cannot recover from the complex $W(X)$ the full mixed Hodge
structure on the rational cohomology of $X$, but only the graded
quotients of its weight filtration ; for this reason we call $W(X)$
the {\it weight complex} of $X$.  We indicate in Section 3.2.4 how
$[X]$ relates to some constructions of Voevodsky in the derived
category of mixed motives \cite{V}. Hanamura told us that he had,
independently, used cohomological descent for higher Chow groups
tensored with ${\Bbb Q}$ to associate mixed motives to arbitrary
varieties.  Concerning Grothendieck's work on Serre's question and
virtual motives, see \cite{Gr3} pp.185 and 191.

\medskip
Our construction gives a new proof of the existence of virtual Betti
numbers for complex varieties, which does not use Hodge theory
(Section 3.3.1). Furthermore we prove that, from $E_2$ on, the weight
spectral sequence of singular cohomology with compact supports and
{\it arbitrary} constant coefficients is independent of choices when
defined by means of (hyper)envelopes ( Theorem 3); for rational
coefficients and arbitrary proper hypercoverings this is a result of
Deligne \cite{D}.  In particular, we obtain a canonical weight
filtration on singular cohomology with compact supports and arbitrary
coefficients.

\medskip
An example of descent comes from the standard square diagram of
varieties associated to a monoidal transform. When the center is
regularly embedded, one can show directly that the corresponding Chow
groups and higher $K$-groups fit into short exact sequences (Theorem
5). This gives an alternative approach to some of the results in the
paper.

\medskip

If linear equivalence is replaced by any other adequate relation on
algebraic cycles (for instance homological or numerical equivalence),
the same results are clearly also valid. Furthermore, motives can also
be defined by taking for correspondences from $X$ to $Y$ the
Grothendieck group $K_0 (X \times Y)$ (when $X$ and $Y$ are smooth and
projective). For any variety $X$ over a field $k$ of characteristic
zero, we get a complex $KW(X)$ of such motives, well-defined up to
homotopy (4.3).  We then define the {\it higher $K$-groups with
compact support} of $X$ (Theorem 7). These groups $K_m^c (X)$, $m\geq
-\dim (X)$, which coincide with the usual $K$-theory groups for
projective non-singular $X$, admit a weight filtration and a pairing
with the $K'$-theory of $X$ (Proposition 8).

\medskip

The paper is organized as follows. In Section 1 we introduce the
Gersten complexes, our version of Manin's identity principle (Theorem
1), motives and descent. In Section 2 we prove the existence of $W(X)$
and derive its main properties: functoriality, Mayer-Vietoris property
and mulitiplicativity (Theorem 2). In Section 3 we use $W(X)$ to
define $[X]$ as well as other invariants of varieties, we define the
weight filtration on cohomology with compact supports (Theorem 3), and
we give some examples.  In Section 4 (which is almost entirely
independent of the previous ones) we study the behaviour of Chow and
$K$-groups under blow up. In Section 5 we define the $K$-theory
spectrum with compact support of varieties (Theorem 7). The Appendix
contains some facts about (co)homotopy limits and pairings which are
needed in Section 5.

\medskip

The results of this paper were presented at the International
Conference on Algebraic K-theory in Paris, July 1994. A new
construction of $W(X)$ has been given by Guillen and Navarro in
\cite{G-N}, where it is denoted $h_c(X)$. They construct also a
canonical complex of motives $h(X)$, corresponding to cohomology
without support.

\bigskip

We thank S. Bloch, L. Illusie, N.Katz, Y. Manin, W. Messing, V.
Navarro, T.  Scholl, J.P. Serre and B. Totaro for helpful discussions.

\bigskip

\heading 1. Motivic Descent \endheading

\subheading{1.1 The Gersten complex}

\medskip

Let $X$ be a scheme of type over a field and $q\geq 0$ an integer. The
{\it Gersten complex} is a chain complex of abelian groups
$R_{q,*}(X)$ such that $$R_{q,i}(X)=\bigoplus_{x\in X,
{\dim}\overline{\{x\}}=q+i}K_{i}({\bold k}(x)),$$ where $K_{i}({\bold
k}(x))$ is the i-th higher K-group of the residue field at the point
$x$ of $X$, and $\overline{\{x\}}$ denotes the Zariski closure of the
set $\{x\}$.  When $q$ varies, these complexes form the $E_1$-term of
the spectral sequence associated to the $K$-theory of the category of
coherent sheaves on $X$, filtered by dimension of support \cite{Q} \S
7,Theorem 5.4 and
\cite{G1}, and the differential is the corresponding $d_1$.
  It is proved in \cite{G1} Th. 7.22 that, if $f:X\to Y$
is a proper map of schemes, then we have
a map of chain complexes: $$
\align f_*:R_{q,*}(X) \to & R_{q,*}(Y)\\ \bigoplus_x
 f_x : \bigoplus_{x\in X, {\dim}\overline{\{x\}}=q+i}K_{i}({\bold
k}(x)) \to & \bigoplus_{y\in Y,
{\dim}\overline{\{y\}}=q+i}K_{i}({\bold k}(y)), \endalign $$ where
$f_x=0$ if $\dim\overline{\{f(x)\}}\}<
\dim\overline{\{x\}}$, and $f_x$ is the norm map associated to the
finite field extension ${\bold k}(x)/{\bold k}(y)$ if
$\dim\overline{\{f(x)\}}\}= \dim\overline{\{x\}}$.  The 0-th homology
group of the complex $R_{q,*}$ is the Chow homology group of dimension
$q$ cycles, $CH_q(X)$, while the other homology groups are the higher
$K$-theory type Chow homology groups $CH_{q,p}(X)=H_p(R_{q,*}(X))$,
which are related to $K_{q-p}(X)$. As shown in \cite{G1} p.276, the
direct sum of these groups for all $p$ and $q$ is a graded module over
the ring $\bigoplus_{p,q}H^p(X,{\Cal K}_q)$, where ${\Cal K}_q$ is the
Zariski sheaf associated to the presheaf $K_q$. In particular, for a
non-singular $X$ of finite type over a field, it is a graded module
over the Chow ring.  \bigskip

\subheading{1.2 Universal acyclicity} \medskip Throughout this paper,
given a field $k$, a variety over $k$ means a reduced scheme of finite
type over ${\text {Spec}}(k)$, and we denote by ${\bold V}$ either the
category of smooth proper varieties over $k$ or the category of smooth
projective varieties over $k$: all statements will be valid for both
definitions of ${\bold V}$ (except in section 5 where we shall require
that varieties in ${\bold V}$ are projective).  Let $\Bbb Z\bold V$ be
the category with the same objects as $\bold V$, but with
$\text{Hom}_{\Bbb Z\bold V}(X,Y)$ equal to the free abelian group
${\Bbb Z}\text{Hom}_{\bold V}(X,Y)$ on $\text{Hom}_{\bold V}(X,Y)$.
For varieties $X$, $Y$ and $Z$, the composition pairing $$
\text{Hom}_{\Bbb Z\bold V}(X,Y)
\times\text{Hom}_{\Bbb Z\bold V}(Y,Z) \to \text{Hom}_{\Bbb Z\bold
V}(X,Z)$$ is bilinear and induced by the usual composition of
morphisms in ${\bold V}$.  Notice that, for each $q\geq 0$, $ X\mapsto
R_{q,*}(X)$ is a covariant functor from $\bold V$ to the category of
chain complexes of abelian groups, and therefore factors through
${\Bbb Z}\bold V$.

Let $\bold C$ be the category of correspondences in $\bold V$, having
the same objects as $\bold V$, but with $\text{Hom}_{\bold
C}(X,Y)=\bigoplus_{i\in I}CH^{{\dim}(Y_i)}(X\times Y_i)$, where
$Y_i,i\in I$, are the connected components of $Y$.  The composition
law is defined as follows: $$\text{Hom}_{\bold C}(X,Y)\times
\text{Hom}_{\bold C}(Y,Z)\to \text{Hom}_{\bold C}(X,Z)$$ $$(\alpha
,\beta)\mapsto
\pi_{XZ*}(\pi^*_{XY}(\alpha)\pi^*_{YZ}(\beta))$$ where $\pi_{XZ}$,
$\pi_{XY}$ and $\pi_{YZ}$ are the projections from $X\times Y\times Z$
to $X\times Z$, $X\times Y$ and $Y\times Z$ respectively.  This
composition is bi-additive, and for any $X\in Ob \bold V$ the class
$1_X =[\Delta_X]$ of the diagonal is the identity in the ring
$\text{End}_{\bold C}(X)$.  There is a {\it covariant} functor
$$\Gamma: {\bold V} \to {\bold C}$$ mapping a morphism $f$ to the
class $ [\Gamma_f]\in
\text{Hom}_{\bold C}(X,Y)$ of the graph of $f$.  This functor factors
through ${\Bbb Z}{\bold V}$. The Chow homology type functors
$CH_{p,q}$ all factor through $\bold C$, with a correspondence acting
on a class $x$ by $\alpha_*:x\mapsto \pi_{Y*}(\alpha \pi^*(x))$
(\cite{So1} Section 1.2).

Recall from \cite{M} p.448 the following lemma:

\proclaim{Lemma 1} Let $f:X\to Y$ be a map in $\bold V$, and $\theta
\in CH^*(Z\times X)$ a correspondence. Then $$(1_Z\times f)_*\theta =
[f]\circ \theta \in CH^{{\dim}(Y)}(Z\times Y)$$ \endproclaim

\demo{ Proof}
 $$\align [f]\circ \theta =
&\pi_{ZY*}(\pi^*_{ZX}(\theta)\pi^*_{XY}([f]))\\ =&
\pi_{ZY*}(\pi^*_{ZX}(\theta)j_*([Z\times X])\\ =&
\pi_{ZY*}(j_*j^*\pi_{Z\times X}(\theta)))\\ =& \pi_{ZY*}(j_*(\theta
))\\ =& (\pi_{ZY}\circ j)_*(\theta ))\\ =& (1_Z\times f)_*(\theta
)\endalign $$ Here $j:Z\times X\to Z\times X\times Y$ is the graph of
$f\circ \pi_X : Z\times X\to Y$.  \qed\enddemo

\proclaim{Theorem 1} Suppose that $$\ldots\to X_2 \overset \delta_2 \to
\to X_1 \overset \delta_1 \to \to X_0$$ is a chain complex in ${\Bbb
Z}{\bold V}$, such that for all $V\in \bold V$ and all $q\geq 0$, the
total complex of the double complex $$\ldots\to R_{q,*}(V\times
X_2)\to R_{q,*}(V\times X_1)\to R_{q,*}(V\times X_0)$$ is acyclic.
Then the complex $$\ldots\to X_2 \overset [\delta_2] \to \to X_1
\overset [\delta_1] \to \to X_0$$ in $\bold C$ has a contracting
homotopy.  \endproclaim

\demo{ Proof } We proceed by induction on $n\geq 0$ to construct
correspondences $h_n:X_n\to X_{n+1}$ such that
$h_{n-1}\delta_n+\delta_{n+1}h_n=1_{X_n}$.

We start with the case of $n=0$. Consider the chain complex $$
\ldots\to X_0\times X_2 \overset \delta_2 \to \to X_0\times X_1
\overset \delta_1 \to \to X_0\times X_0 $$ and write $d= {\dim}(X_0)$.
By hypothesis, the total complex associated to the double complex $$
\ldots R_{d,*}(X_0\times X_2) \overset \delta_2 \to \to
R_{d,*}(X_0\times X_1 ) \overset \delta_1 \to \to R_{d,*}(X_0\times
X_0) $$ is acylic. In particular $$ R_{d,1}(X_0\times X_0 )\oplus
R_{d,0}(X_0\times X_1 )\to R_{d,0}(X_0\times X_0 ) $$ is surjective.
Hence $$ (1_{X_0}\times \delta_1)_*:R_{d,0}(X_0\times X_1 )\to
CH_d(X_0\times X_0)=CH^d(X_0\times X_0) $$ is surjective. In
particular there exists a cycle $\eta_0$ on $X_0\times X_1$ such that
$$ (1_{X_0}\times \delta_1)_*(\eta_0)=[1_{X_0}](=[\Delta_{X_0}]). $$
By Lemma 1, if we set $h_0=[\eta_0]$, the correspondence represented
by $\eta_0$, then $$ \delta_1\circ h_0 = 1_{X_0}\quad .  $$

Suppose now that $n \geq 1$ and that correspondences $h_i:X_i\to
X_{i+1}$ for $i=0\ldots n-1$ have been constructed, such that $$
h_{i-1}\circ \delta_i +
\delta_{i+1}\circ h_i = 1_{X_i}\quad .  $$
 Since, for all $p, q\geq 0$, $CH_{q,p}$ factors
 through $\bold C$, it follows that for all varieties $ V$ we have
maps: $$ (1_V\times
\delta_i)_*:CH_{q,p}(V\times X_i )\to CH_{q,p}(V\times X_{i-1} ) $$
for all $i$, and maps $$ (1_V\times h_i)_*:CH_{q,p}(V\times X_i )\to
CH_{q,p}(V\times X_{i+1} ) $$ for $i\leq n-1$, such that, for $i\leq
n-1$, $$ (1_V\times h_{i-1})_*\circ (1_V\times \delta_i)_* +
(1_V\times \delta_{i+1})_*\circ (1_V\times h_i)_* =
1_{CH_{q,p}(V\times X_i)}\quad .  $$ It follows that, for $i\leq n-1$,
$$ H_i(*\mapsto CH_{q,p}(V\times X_* ))=0 $$ Now by hypothesis, for
any $q\geq 0$, the total complex $\int R_{q,*}(V\times X_*)$
associated to the double complex $R_{q,*}(V\times X_*)$ is acyclic.
Consider the spectral sequence $$E^2_{m,p}=H_m(*\mapsto
CH_{q,p}(V\times X_*))\Rightarrow H_{m+p}(\int R_{q,*}(V\times
X_*)).$$ From the induction hypothesis and the fact that the functors
$CH_{q,p}$ factor through $\bold C$, we deduce that $E^2_{m,p}=0$ when
$0\leq m \leq n-1$ and $p\geq 0$.  Therefore $E^r_{m,p}=0$ whenever
$r\geq 2$, $0\leq m \leq n-1$ and $p\geq 0$, and the only part of
$\bigoplus_{m+p=n}E^2_{m,p}$ which can possibly be non-zero is
$E^2_{n,0}=H_n(*\mapsto CH_q(V\times X_*))$. Using the induction
hypothesis again, we have that all the differentials out of
$E^2_{n,0}$ are zero. On the other hand, since $p=0$, there are no
differentials {\it into} $E^2_{n,0}$. Hence $$
E^2_{n,0}=E^\infty_{n,0}=0\quad , $$ since $H_n(\int R_{q,*}(V\times
X_*)=0$.  Therefore, for all varieties $V$, and all $q\geq 0$, we have
that $$ CH_q(V\times X_{n-1})\leftarrow CH_q(V\times X_{n})\leftarrow
CH_q(V\times X_{n+1})$$ is exact in the middle.

Now take $V=X_n$ and $q= {\dim}(X_n)$. Consider the element $$
1_{X_n}-h_{n-1}\circ \delta_n \in CH_q(X_{n}\times X_{n})\simeq
CH^q(X_{n}\times X_{n})\quad , $$ By Lemma 1, the image of this class
under $$(1_{X_n}\times \delta_n)_*:CH_q(X_{n}\times X_{n})\to
CH_q(X_{n}\times X_{n-1})$$ is the correspondence $$
\align
\delta_n\circ(1_{X_n}-h_{n-1}\circ \delta_n) = &\delta_n-\delta_n\circ
h_{n-1}\circ \delta_n \\ = & \delta_n -(1_{X_{n-1}}-h_{n-2}\circ
\delta_{n-1})\circ\delta_n \\ = & 0 \quad \text{since}\quad
\delta_{n-1} \circ \delta_n =0 .\endalign $$ Hence, by the exactness
above, there exists a cycle $$\eta_n\in CH_q(X_{n}\times X_{n+1})$$
such that
$$(1_{X_n}\times\delta_{n+1})_*\eta_n=1_{X_n}-h_{n-1}\circ\delta_n.$$
This cycle represents a correspondence $$h_n:X_n\to X_{n+1}$$ which
satisfies, by Lemma 1, the identity $$\delta_{n+1}\circ
h_n+h_{n-1}\circ \delta_n=1_{X_n}$$ and we are done.  \qed\enddemo

\subheading{1.3 Motives} \medskip Let $\bold M$ be the category of
(pure effective) {\it Chow motives} over $k$. It is obtained from the
category $\bold C$ of correspondences by inverting the arrows and
adding the images of projectors (cf. \cite{M1}, \cite{Kl} ).  An
object of $\bold M$ is a pair $(X,p)$ where the variety $X$ lies in
$\bold V$ (i.e. it is smooth and either projective or proper over $k$)
and $p \in End_{\bold C}(X)$ satisfies $p^2 = p$. A morphism from
$(X,p)$ to $(Y,q)$ is an element in $qHom_{\bold C}(Y,X)p$ (see
\cite{J} or \cite{Sc} for this definition). There is a contravariant
functor from $\bold V$ to $\bold M$ mapping a variety $X$ to the Chow
motive $M(X)=(X,1_X)$ and a morphism $f$ to the transpose of
$[\Gamma_f]$. The category $\bold M$ is pseudo-abelian, i.e. it is
additive and projectors have images. Disjoint union and product of
varieties can be extended to motives where they are denoted $\oplus$
and $\otimes $ respectively.

More generally, given any equivalence relation $\sim$ on algebraic
cycles which is adequate in the sense of \cite{Kl}, one may substitute
to Chow groups the groups of cycles modulo $\sim$ in the above
construction, getting a category of motives $\bold M_{\sim}$. Basic
examples of adequate equivalence relation are linear equivalence,
homological equivalence and numerical equivalence, written $ num$ in
what follows. Notice that for any ${\sim}$ there exist covariant
functors $\bold M \to \bold M_{\sim}$ and $\bold M_{\sim} \to \bold
M_{num}$ (\cite{Kl} Proposition 3.5).

Given any chain complex $X.$ in ${\Bbb Z}{\bold V}$, we shall denote
by $M(X.)$ (resp. $M_{\sim}(X.)$) the corresponding cochain complex in
$\bold M$ (resp. $\bold M_{\sim}$).

\proclaim{Corollary 1} Let $f:X.\to Y. $ be a morphism of chain
complexes in ${\Bbb Z}{\bold V}$ such that for all $V\in {\bold V}$
and all $q\geq 0$, the induced map of double complexes
$$R_{q,*}(V\times X_*)\to R_{q,*}(V\times Y_*) $$ induces a
quasi-isomorphism of total complexes.  Then the corresponding map
$X.\to Y.$ of chain complexes in $\bold C$ is a homotopy equivalence.

Consequently, for any choice of an adequate equivalence relation, the
induced map $f^*:M_{\sim}(Y.)\to M_{\sim}(X.)$ is a homotopy
equivalence.
\endproclaim

\demo{Proof} To prove this corollary we shall use a general theorem of
Verdier \cite{V}.  Let $\bold A$ be an arbitrary additive category.
For any map $f:X.\to Y. $ of complexes in $\bold A$, denote by
$Cone(f). = C(f).$ the mapping cone of $f$. According to \cite{V} II,
Proposition 1.3.2, the category $Hot(\bold A)$ of complexes in $\bold
A$ up to homotopy is a triangulated category, in which the triangles
are the diagrams isomorphic in $Hot(\bold A)$ to diagrams of the form
$$X.\overset f \to \to Y. \to C(f). \to X.[1]$$ where $f:X.\to Y.$ is
any morphism of complexes in $\bold A$.

  It follows from this that if $f:A.\to B.$ is a map between chain
complexes in $\bold A$ the following two statements are equivalent:

\item{$\bullet$} The map is a chain homotopy equivalence (i.e. there
exist a chain map $g:B.\to A.$ and homotopies $g\circ f\sim 1_{A.}$
and $f\circ g \sim 1_{B.}$);

\item{$\bullet$} The mapping cone $C(f).$ is contractible.

Indeed, if $f$ is a homotopy equivalence and if ${\text {id}} : X.
\rightarrow X.$ is the identity map, there exists a morphism of
triangles in ${\text {Hot}}({\bold A})$

$$\CD X.  @>{ \text{id}}>> X.  @>>> C({\text{id}}). @>>> X.[1] \\
@VV{\text{id}}V @VV{f}V@VV{\phi}V @VV{\text{id}}V\\ X.  @>{f}>> Y.
@>>> C(f).  @>>> X.[1]
\endCD $$ (property TRIII in \cite{V} II 1.1.1).
 From Cor. 1.2.5 in \cite{V} II, it follows that
 $\phi$ is an isomorphism in ${\text {Hot}}({\bold A})$.
Since $C({\text {id}}).$ is contractible, it follows that $C(f).$ is
contractible.

\medskip

Conversely, if $C(f).$ is contractible it is isomorphic to 0 in
${\text {Hot}}({\bold A})$ and we get a triangle $X. \overset{f}\to
\rightarrow Y. \rightarrow 0 \rightarrow X. [1]$, hence a triangle $0
\rightarrow X. \overset {f}\to \rightarrow
 Y. \rightarrow 0$ (TRII in \cite{V} II 1.1.1). Therefore $f$ is an
isomorphism in ${\text {Hot}}({\bold A})$, by TRI in \cite{V} II 1.1.1
and the uniqueness of the cone, \cite{V} II Cor. 1.2.6.

To prove the corollary, notice that $R_{q,*}$ commutes with the
formation of mapping cones of complexes in ${\Bbb Z}{\bold V}$. So,
under our assumptions, the complex $C(f).$ in ${\Bbb Z}{\bold V}$
satisfies the hypotheses of Theorem 1.  Therefore $C(f).$ is
contractible when viewed as a complex in $\bold C$, and hence the map
$f:X.\to Y.$ is a homotopy equivalence of complexes in $\bold C$.
\qed\enddemo \bigskip

\subheading{1.4 Envelopes} \subheading{1.4.1 } An {\it envelope}
$p:X\to Y$ is a proper map of schemes such that for all fields $F$,
the induced map $X(F)\to Y(F)$ is surjective. A {\it hyperenvelope} is
a map $p:X.\to Y.$ of simplicial schemes which is proper in each
degree and which is a hypercovering for the Grothendieck topology in
which the covering maps are envelopes. Specifically, for all $i\geq 0$
the map $$ X_i\to (cosk^{Y.}_{i-1}sk_{i-1}(X.))_i$$ is an envelope
(see \cite{G2}).  Equivalently, a hyperenvelope is a proper map of
simplicial schemes such that for each field $F$ the map of simplicial
sets $X.(F)\to Y.(F)$ is a trivial Kan fibration \cite{Ma}.

 Let $f:Y.\to X.$ be a morphism of simplicial schemes. A {\it
hyperenvelope} $h:\tilde{f}\to f$ is a map in the category of arrows
of simplicial schemes, {\rm i.e.} a commutative square $$ \CD
\tilde{Y}.@>\tilde{f}>>\tilde{X}.\\ @Vh_YVV @VVh_XV \\ Y. @>>f> X.
\endCD \quad ,$$ such that $h_Y$ and $h_X$ are hyperenvelopes.

If $k$ is a field of characteristic zero, any projective variety $X$
over $k$ has a {\it non-singular envelope}, i.e.  an envelope lying in
${\bold V}$.  This can be shown by induction on the dimension of $X$:
by Chow's lemma and Hironaka's resolution of singularities one may
find $X'\in {\bold V}$ and a proper map $X'\to X$, which is an
isomorphism over a dense open subset $U\subset X$; the disjoint union
of $X'$ with a non-singular envelope of $X-U$ is a non-singular
envelope of $X$.  From this one gets by induction on the simplicial
degree that, when viewed as a constant simplicial scheme, $X$ admits a
hyperenvelope $Z.$ which is non-singular in the sense that each $Z_i,
i\ge 0 $ lies in ${\bold V}$.

More generally, the following holds:

\proclaim{Lemma 2} If $k$ is a field of characteristic zero any
simplicial variety over $k$ has a non-singular envelope.  Furthermore,
given a morphism of simplicial varieties $f:Y.\to X.$ over $k$, there
exists a non-singular hyperenvelope $h:\tilde{f}\to f$.  \endproclaim

\demo{Proof} To prove the first assertion consider a simplicial
variety $X.$. One defines a non-singular hyperenvelope $f:Y.\to X.$ by
induction on $n\geq 0$, as follows (see \cite{SD} 5.1.3 and \cite{D}
Section 6.3 for the analogous result for cohomological descent).
Start by choosing a non-singular envelope $f_0:Y_0\to X_0$.  Suppose
now that $n>0$, and that the $(n-1)$-skeleton of $Y.$ has been
constructed together with maps $f_i:Y_i\to X_i$ for $i=0\ldots n-1$
which define a hyperenvelope of $(n-1)$-truncated simplicial schemes.
Let $Z_n\to cosk_{n-1}^{X.}(sk_{n-1}Y.)_n$ be a non-singular envelope.
Set $$Y_n = Z_n \coprod \coprod_{j=0}^{n-1}Y_{n-1}\quad .$$ For
$i=0\ldots n-1$, let the degeneracy map $s_i:Y_{n-1}\to Y_n$ be the
identity map onto the $i$-th summand in the coproduct. The face maps
$d_i:Y_n\to Y_{n-1}$ are then defined as follows.  On $Z_n$ they are
the composition of the envelope $Z_n\to
cosk_{n-1}^{X.}(sk_{n-1}(Y.))_n$ with the face maps
$d_i:cosk_{n-1}^{X.}(sk_{n-1}Y.)_n\to Y_{n-1}$.  On the $j$-th summand
in $\coprod_{j=0}^{n-1}Y_{n-1}$, we set $d_i=s_kd_l$ where $k=j\pm 1$
or $k=j$ and $l=i-1$ or $l=i$, according to the formula in definition
1.1 (iii) of \cite{Ma} .  Let $f_n:Y_n\to X_n$ be the map which on
$Z_n$ is the composition of the envelope $Z_n\to
cosk_{n-1}^{X.}(sk_{n-1}Y.)_n$ with the natural map
$cosk_{n-1}^{X.}(sk_{n-1}Y.)_n\to X_n$.  On the $i$-th summand of the
coproduct $Z_n\to cosk_{n-1}^{X.}(sk_{n-1}Y.)_n$ we set
$f_n=s_jf_{n-1}$. We define the $n$-skeleton of $Y.$ to be the
$n$-truncated simplicial scheme obtained by this process, and the
sequence of maps $f_i$ for $i=0\ldots n$ are a hyperenvelope
$sk_n(Y.)\to sk_n(X.)$.

To prove the second assertion, let $\tilde {X.}$ be a non-singular
hyperenvelope of $X.$.  Now choose $\tilde{Y.}$ to be a non-singular
hyperenvelope of the pull- back $Y.\times_{X.}\tilde{X.}$, with
$\tilde{f}:\tilde{Y.}\to
\tilde{X.}$ the obvious composition.  We get a non-singular envelope
$h:\tilde{f}\to f$, i.e. a commutative square $$ \CD
\tilde{Y}.@>\tilde{f}>>\tilde{X}.\\ @Vh_YVV @VVh_XV \\ Y. @>>f> X.
\endCD \quad $$ where $h_X$ and $h_Y$ are hyper-envelopes and, for all
$i\ge 0$, $\tilde{Y}_i$ and $\tilde{X}_i$ lie in ${\bold V}$.
\qed\enddemo
\bigskip

The following are some basic properties of hyperenvelopes, which
follow in a straightforward fashion from the definitions :

\item{$\bullet$} The composition of two hyperenvelopes is again a
hyperenvelope.

\item{$\bullet$} If $f:X.\to Y.$ is a hyperenvelope, and $g:Z.\to X.$
is a map of simplicial schemes, then the pull-back of $f$ along $g$ is
again a hyper-envelope. The same is true for envelopes of arrows of
simplicial schemes.

\item{$\bullet$} The fibre product of two hyperenvelopes is again a
hyperenvelope (this follows from the previous two assertions).

\subheading{1.4.2} If $X.$ is a simplicial scheme with proper face
maps, then for each $q\geq 0$ we have the Gersten complexes
$R_{q,*}(X.)$, by which we mean the total complex associated to the
double complex $(i,j)\mapsto R_{q,i}(X_j)$.  Given a proper map
$f:X.\to Y.$ of simplicial schemes (with proper face maps), there is a
push-forward map $f_*:R_{q,*}(X.)\to R_{q,*}(Y.)$, making
$R_{q,*}(X.)$ a covariant functor from the category of simplicial
schemes with proper face maps to the category of chain complexes.
Envelopes then have the following descent property.

\proclaim{Proposition 1} If $f:X.\to Y.$ is a hyperenvelope of
simplicial schemes with proper face maps, then $f_*:R_{q,*}(X.)\to
R_{q,*}(Y.)$ is a quasi-isomorphism.  \endproclaim \demo{Proof} For
the case when $Y.$ is the constant simplicial scheme associated to a
scheme $Y$ ( which is the only case that we need in this paper) this
is proved in \cite{G2} Theorem 4.3.  The general case follows by an
extension of the argument of op. cit..  Specifically, if we let $X[i].
= cosk^{Y.}_{i}(sk_{i}(X.))$ when $i\geq 0$ and $X[-1]. = Y.$, one can
show that the map $X[i].
\rightarrow X[i-1].$ induces a quasi-isomorphism of Gersten complexes
for all $i\geq 0$: the argument in \cite{SD}, proof of Theorem 3.3.3,
reduces this assertion to the descent theorem in the case $Y.$ is
constant and $X. = cosk_0^{Y_0} (X_0)$, i.e. to Step II in
\cite{G2}
 Theorem 4.3. \qed \enddemo

\bigskip

\subheading {1.4.3 } Given a simplicial object $X.$ in $\bold V$,
write ${\Bbb Z}X.$ for the complex $$ \ldots X_i\overset
\sum_{j=0}^{i} (-1)^j d_j \to \to X_{i-1}\ldots $$ which we view as a
complex in ${\Bbb Z}\bold V$, and $M(X.)$ for the corresponding
complex of motives.

\proclaim{Proposition 2} Let $p:X.\to Y.$ be a map of simplicial
objects in $\bold V$ which is a hyperenvelope. Then the associated map
of complexes in $\bold M$: $$p^*:M(Y.)\to M(X.)$$ is a homotopy
equivalence.
\endproclaim

\demo{Proof} Since $p$ is a hyper-envelope, for any variety $V$ so is
$1_V\times p:V\times X.\to V\times Y.$ , see 1.4.1. or \cite{G2},
Lemma 3.1 i).  Therefore, by Proposition 1, for any $q\geq 0$ the map
of Gersten complexes $$p_*: R_{q,*}(V\times X.)\to R_{q,*}(V\times
Y.)$$ is a quasi-isomorphism.  Consider the map of chain complexes
${\Bbb Z}X.\to {\Bbb Z}Y.$ induced by $p$.  Since ${\Bbb Z}(V\times
X.)  = V\times {\Bbb Z}X.$ and the formation of $ R_{q,*}(X.)$ can be
done using either $X.$ or ${\Bbb Z}X.$, we conclude from Corollary 1
that $p^*:M(Y.)\to M(X.)$ is a homotopy equivalence. \qed \enddemo

\bigskip

\heading 2. The weight complex of arbitrary varieties \endheading

\subheading{2.1} Let $k$ be a field of characteristic zero and
 $X$ a (possibly singular) variety over $k$. Choose a {\it
compactification} $\bar{X}$ of $X$, by which we mean any complete
variety containing $X$ as an open subvariety, with complement
$Y=\bar{X}-X$. Then for any non-singular hyper-envelope
$p:\tilde{X}.\to\bar{X}$ there exists a non-singular hyper-envelope
$q:\tilde{Y}.\to Y$, and a map $\tilde{j}:\tilde{Y}.\to\tilde{X}.$
such that we have a commutative diagram, where $j:Y\to \bar{X}$ is the
inclusion: $$ \CD
\tilde{Y}.@>\tilde{j}>> \tilde{X}.\\
    @VqVV @VVpV \\ Y @>>j> \bar{X}\endCD $$ We associate to this data
the cochain complex $S\left(\tilde{j}\right) =
Cone(\tilde{j}^*)^.[-1]$, where $\tilde{j}^*$ is the induced map of
cochain complexes $M(\tilde{X}.)\to M(\tilde{Y}.)$; in particular, for
all $n\ge 0$, $S\left(\tilde{j}\right)^n= M(\tilde{X}_n)\oplus
M(\tilde{Y}_{n-1})$.  Let $W(X)$ be the class of
$S\left(\tilde{j}\right)$ in the homotopy category ${Hot}(\bold M)$ of
cochain complexes in $\bold M$.

\proclaim{Theorem 2} In ${Hot}(\bold M)$
 the complex $W(X)$ is independent, up to canonical isomorphism, of
the choices made. It has the following properties: \smallskip
\item{i)} $W(X)$ is isomorphic to a bounded complex
of the form $$M(X_0) \rightarrow M(X_1) \rightarrow \cdots \rightarrow
M(X_k)$$ where, for all $i$, $0\leq i\leq k$, $X_i$ is a variety lying
in $\bold V$ and ${\dim} (X_i) \leq \dim(X)-i$; in particular $k\leq
\dim(X)$.

\smallskip

\item{ii)} Any proper map $f:X\rightarrow X'$ of varieties induces a
morphism $f^* :W(X')\rightarrow W(X)$ in $Hot (\bold M)$. Given two
composable proper maps $f$ and $g$, on has $(fg)^* =g^* f^*$.  Any
open immersion $i:X\rightarrow X'$ induces a morphism $i_*
:W(X)\rightarrow W(X')$ in $Hot (\bold M)$ and, given two composable
open immersions $i$ and $j$, $(ij)_* =i_* j_*$.

\smallskip

\item{iii)} Let $X$ be a variety, and suppose that $i:U\rightarrow X$
is an open immersion, with complement the closed immersion
$f:T\rightarrow X$.  Then there is a canonical triangle in
${Hot}(\bold M)$: $$ W(U)\overset {i_*}
\to \to W(X)\overset {f^*} \to
\to W(T) \to W(U)[1]\quad .  $$

\smallskip

\item{iv)} Assume $X$ is the union of two closed subvarieties $A$ and
$B$. Then there is a canonical triangle in ${\text {Hot}} ({\bold M})$
$$ W(X) \rightarrow W(A) \oplus W(B) \rightarrow W(A \cap B)
\rightarrow W(X)[1].  $$ If $X$ is the union of two open subvarieties
$U$ and $V$, there is a canonical triangle $$ W(X)[-1] \rightarrow W(U
\cap V) \rightarrow W(U) \oplus W(V) \rightarrow W(X).  $$ \smallskip

\item{v)} If $X$ and $Y$ are quasi-projective varieties, then
$$W(X\times Y)=W(X) \otimes W(Y).$$ \endproclaim

\bigskip The complex $W(X)$ in $Hot({\bold M})$ will be called the
{\it weight complex} of the variety $X$. As the proof in sections 2.2
to 2.6 will show, Theorem 2 is already true in the homotopy category
of bounded complexes in ${\bold C}$ instead of $Hot({\bold M})$.

\bigskip \subheading{2.2}
 Let ${\bold P}$ be the category of proper varieties over $k$ and
${\text {Ar}} ({\bold P})$ the category of arrows in ${\bold P}$,
where a morphism $g : f_1
\rightarrow f_2$ is defined to be a commutative square
 $$\CD Y_1 @>{f_1}>> Z_1 \\ @Vg_YVV @VVg_ZV
 \\ Y_2 @>>{f_2}> Z_2 \endCD \leqno (2.1) $$ in ${\bold
P}$. As a preliminary to the proof of Theorem 2, we shall define a
functor $T$ from ${\text {Ar}} ({\bold P})$ to ${\text {Hot}} ({\bold
M})$.

\medskip

First suppose that a morphism $g : f_1 \rightarrow f_2$ as above is
given and that we have two non-singular hyperenvelopes $h_i :
\widetilde{f}_i \rightarrow f_i$, $i=1,2$, i.e. commutative squares $$
\CD \tilde{Y}_{.,i}@>\tilde{f}_i>> \tilde{Z}_{.,i}\\
    @Vh_{Y,i}VV @VVh_{Z,i}V \\ Y_i @>f_i>> Z_i \endCD $$ in which
$h_{Y,i}$ and $h_{Z,i}$ are non-singular hyperenvelopes. Consider the
fiber product $\widetilde{f}_1 \times_{f_2} \widetilde{f}_2$ and the
projections $p_i :
\widetilde{f}_1 \times_{f_2} \widetilde{f}_2 \rightarrow
\widetilde{f}_i$, $i=1,2$. By 1.4.1, $p_1$ is an hyperenvelope;
however it may not be non-singular. So let $\pi : \widetilde f
\rightarrow \widetilde{f}_1
 \times_{f_2} \widetilde{f}_2$ be a non-singular hyperenvelope. If
$a=p_1 \circ \pi$ and $b=p_2 \circ \pi$ we get a commutative diagram
$$ \CD \tilde{f}@>b>> \tilde{f}_2\\ @VaVV @VVV \\
\tilde{f_1}@.@VVh_2V\\ @Vh_1VV @VVV \\ {f}_1 @>g>> {f}_2 .  \endCD
\leqno (2.2)$$ Since $a$ is an hyperenvelope we deduce from
Proposition that $S(a) : S\left( \widetilde{f}_1 \right)
\rightarrow S \left( \widetilde f \right)$ is an isomorphism in
${\text {Hot}} ({\bold M})$. Therefore we get a map
$\theta_{\widetilde f} = S(a)^{-1} \ S(b)$ from $S\left(
\widetilde{f}_2 \right)$ to $S \left( \widetilde{f}_1 \right)$.

\medskip

This map does not depend on the choice of $\pi$. Indeed, assume $$\pi'
: \widetilde{f}' \rightarrow \widetilde{f}_1 \times_{f_2}
\widetilde{f}_2$$
 is another non-singular hyperenvelope. Choose a non-singular
hyperenvelope $\widetilde{f}'' \rightarrow \widetilde f \times_{(
\widetilde{f}_1 \times_{f_2} \widetilde{f}_2 )}^{} \widetilde{f}'$. We
have commutative squares

$$ \CD \tilde{f''}@>k'>> \tilde{f'}\\ @VkVV @VVa'V \\ \tilde{f} @>a>>
\tilde{f}_1 \endCD $$ and $$ \CD
\tilde{f''}@>k'>> \tilde{f'}\\
    @VkVV @VVb'V \\ \tilde{f} @>b>> \tilde{f}_2 \endCD $$ where $k$,
$k'$, $a$ and $a'$ are hyperenvelopes. It follows that $$
\eqalign{ \theta_{\widetilde{f}'} = & \ S(a')^{-1} \ S(b') = S(a)^{-1}
\ S(k)^{-1} \ S(k') \ S(b') \cr = & \ S(a)^{-1} \ S(b) =
\theta_{\widetilde f} \ . \cr } $$ Thus we have a canonical map $$
T(g) : S\left( \widetilde{f}_2 \right) \rightarrow S \left(
\widetilde{f}_1 \right) .  $$ When $g$ is the identity map, $p_2$ is
also an hyperenvelope, therefore $S(b)$ and $T(g)$ are isomorphisms in
${\text {Hot}} ({\bold M})$.

\medskip

We shall now check that $T(g)$ is compatible with composition.
Consider morphisms $g:f_1 \rightarrow f_2$ and $k:f_2 \rightarrow f_3$
in ${\text {Ar}} ({\bold P})$ as well as non-singular hyperenvelopes
$h_i : \widetilde{f}_i \rightarrow f_i$, $i=1,2,3$. From the argument
above, we get a commutative diagram $$ \CD
\tilde{f''}@>{\beta}>>\tilde{f}@>b>> \tilde{f}_3\\ @V{\alpha}VV@VVaV
@VVV\\ \tilde{f'}@>{b'}>> \tilde{f}_2@.@VVV\\ @V{a'}VV @VVV @VVV\\
\tilde{f_1}@.@VVV@VVV\\ @VVV @VVV @VVV\\ f_1 @>g>> f_2 @>k>> f_3\\
\endCD\leqno (2.3) $$ where all vertical maps are hyperenvelopes,
$\widetilde{f}''$ being defined as a non-singular hyperenvelope of
$\widetilde{f}'
\times_{\widetilde{f}_2} \widetilde{f}$; since the map $a$ is a
hyperenvelope, the same is true of $\alpha:\widetilde{f}'' \rightarrow
\widetilde{f}'$ . The obvious morphism $\widetilde{f}'' \rightarrow
\widetilde{f}_1 \times_{f_3} \widetilde{f}_3 = \left( \widetilde{f}_1
\times_{f_2} \widetilde{f}_2 \right) \times_{\widetilde{f}_2} \left(
\widetilde{f}_2 \times_{f_3} \widetilde{f}_3 \right)$ is the
composition of a hyperenvelope and the fiber product of two
hyperenvelopes and is therefore also a hyperenvelope. So, from the
commutativity of (2.3) we deduce that $$T(kg) =
\theta_{\widetilde{f}''} = S( a'\alpha)^{-1}S(b\beta) =S(
a')^{-1}S(b')S( a)^{-1}S( b) =\theta_{\widetilde{f}'} \
\theta_{\widetilde{f}} = T(k) \ T(g).$$

\medskip

We conclude from this discussion that, once we fix a choice ${\bold
h}=(h_f)$ of a non-singular hyperenvelope $h_f : \widetilde f
\rightarrow f$ for each arrow $f : Y \rightarrow Z$ in ${\text {Ar}}
({\bold P})$, there is a contravariant functor $$ T_{\bold h} : {\text
{Ar}} ({\bold P}) \rightarrow {\text {Hot}} ({\bold M}) $$ with
$T_{\bold h} (f) = S
\left( \widetilde f \right)$ and $T_{\bold h}
\left( g : f_1 \rightarrow f_2 \right) = T(g)$. Given two different
choices ${\bold h}=(h_f)$ and ${\bold h}' = (h'_f)$, the maps $T(1_f)
: S \left( \widetilde f \right) \rightarrow S \left( \widetilde{f}'
\right)$ define an isomorphism of functors $T_{\bold h} \overset
{\sim} \to \rightarrow T_{{\bold h}'}$.  We will therefore suppress
the choice ${\bold h}$ from the notation and just write $T$ for
$T_{\bold h}$ for some fixed ${\bold h}$.

\medskip

The functor $T$ has the following property. Consider a morphism $g:
f_1 \rightarrow f_2$ in ${\text {Ar}} ({\bold P})$, i.e. a commutative
square like (2.1).  We say that $g$ is {\it Gersten acyclic} when the
following property holds: if $C.$ is the complex of varieties in
${\bold P}$ $$ Y_1 \overset {(g_Y ,f_1)}
\to \rightarrow Y_2 \oplus Z_1 \overset {f_2 -g_Z}
 \to \rightarrow Z_2 , $$ for any integer $q\geq 0$
and any variety $V$ in ${\bold V}$ the total complex of
corresponding Gersten complexes $R_q (V \times C.)$ is acyclic. Now we
claim that, when $g$ is Gersten acyclic, the map $T(g) : T(f_2)
\rightarrow T(f_1)$ is an isomorphism in ${\text {Hot}} ({\bold M})$.
Indeed, if we consider the product of the diagram (2.2) above by the
identity map $1_V$ on $V$, all morphisms $1_V \times a$, $1_V \times
h_1$, $1_V \times g$ and $1_V \times h_2$ induce quasi-isomorphisms of
Gersten complexes, so the same is true for $1_V \times b$, hence, by
Theorem 1, $S(b)$ is an isomorphism in ${\text {Hot}} ({\bold M})$ and
$T(g) = S(a)^{-1} \ S(b)$ is also an isomorphism.  \bigskip

\subheading{2.3} Given a variety $X$ and a complete variety $\overline
X$ containing $X$ as a Zariski open set, if we write $j: Y \rightarrow
\overline X$ for the inclusion of the complement of $X$, we define $$
W(X) = T(j) = S \left( \widetilde j \right) , $$ for $\widetilde j$
any non-singular hyperenvelope of $j$. We shall prove that $W(X)$ is,
up to canonical isomorphism, independent of choices and contravariant
for proper morphisms.

\medskip

So let $f:X_1 \rightarrow X_2$ be a proper morphism of varieties, each
$X_i$ being equipped with a compactification $X_i \hookrightarrow
\overline{X}_i$ with complement $j_i : Y_i = \overline{X}_i - X_i
\rightarrow \overline{X}_i$.  Consider the Zariski closure
$\overline{X}_f$ of the graph of $f$ in $\overline{X}_1 \times
\overline{X}_2$ and $j_f : \overline{X}_f - X \rightarrow
\overline{X}_f$ the inclusion of the complement. The projections
$\overline{X}_1 \times \overline{X}_2 \rightarrow \overline{X}_i$
induce maps $\pi_i : j_f \rightarrow j_i$. Since the map
$\overline{X}_f \rightarrow
\overline{X}_1$ induces an isomorphism $\overline{X}_f - Y_f
 \rightarrow \overline{X}_1 -Y_1$, the morphism
 $\pi_1 : j_f \rightarrow j_1$ in ${\text {Ar}}
({\bold P})$ is Gersten acyclic. Indeed, for all $q\geq 0$ and $V$ in
${\bold V}$, there is a commutative diagram of Gersten complexes $$
\CD 0@>>> R_{q,*}(V\times Y_f)@>>> R_{q,*}(V\times\overline{X}_f )@>>>
R_{q,*}(V\times(\overline{X}_f - Y_f))@>>> 0\\ @.@VVV @VVV @VV||V\\
0@>>> R_{q,*}(V\times Y_1)@>>> R_{q,*}(V\times\overline{X}_1) @>>>
R_{q,*}(V\times(\overline{X}_1 - Y_1))@>>> 0,\\ \endCD $$ where the
exactness of the rows follows easily from the definition of Gersten
complexes: the horizontal maps are push-forward for closed immersions
and pull-back for open immersions, the exactness degreewise comes from
the fact that, given a scheme and a closed subset, any point lies
either in the closed subset or in its open complement.  It follows
that $T(\pi_1)$ is an isomorphism. We now define $$ W(f) : W(X_2) =
T(j_2) \rightarrow W(X_1) = T(j_1) $$ to be $T(\pi_1)^{-1} \ T(\pi_2)$
(note that $W(f)$ depends on the choice of $j_1$ and $j_2$). When $f$
is an isomorphism, $\pi_2 : j_f \rightarrow j_2$ is also Gersten
acyclic and $W(f)$ is an isomorphism.

\medskip

To check that $W(f)$ is compatible with composition, consider proper
morphisms $f:X_1 \rightarrow X_2$ and $g:X_2 \rightarrow X_3$,
together with compactifications $X_i \rightarrow \overline{X}_i$,
$i=1,2,3$, of complements $j_i : \overline{X}_i - X_i \rightarrow
\overline{X}_i$.

\medskip

Let $\overline{X}_{(f,g)}$ be the Zariski closure of the image of
$X_1$ under the map $(1,f,gf)$ into $\overline{X}_1 \times
\overline{X}_2 \times \overline{X}_3$ and $j_{(f,g)} :
\overline{X}_{(f,g)} - X_1 \rightarrow \overline{X}_{(f,g)}$ its
complement. Then we have a commutative diagram in ${\text {Ar}}
({\bold P})$ $$ \CD j_{(f,g)}@>b>> j_g@>{\pi}_3>> j_3\\ @VaVV
@VV{\pi}'_2V @.\\ j_f@>{\pi}_2>> j_2@.\\ @V{\pi}_1VV @. @.\\ j_1@.@.\\
\endCD $$ induced by the obvious projections. All vertical maps in
this diagram are Gersten acyclic, hence $T$ turn them into
isomorphisms.  Furthermore the composite maps $j_{(f,g)} \rightarrow
j_1$ and $j_{(f,g)}
\rightarrow j_3$ in this diagram factor though the projection
 $j_{(f,g)} \rightarrow j_{gf}$, which is
also Gersten acyclic.  Therefore $u_1: j_{gf}
\rightarrow j_1$ is Gersten acyclic and we have
 a map $u_3: j_{gf} \rightarrow j_3$.  Now we compute $$ W(f) \ W(g) =
(T({\pi}_1)^{-1}T({\pi}_2))(T({\pi}'_2)^{-1}T({\pi}_3))$$
$$=T({\pi}_1a)^{-1}T({\pi}_3b)=T(u_1)^{-1} T(u_3)=W(gf) $$ (from
$W(j_3)$ to $W(j_1)$).

\medskip

Thus, once we fix a choice $\bold {j}=(j_X)$ of compactifications $X
\hookrightarrow \overline X$ with complement $j_X : \overline X - X
\hookrightarrow \overline X$ for all varieties $X$ of finite type over
$k$, we get a contravariant functor $W_{\bold {j}}$ from varieties and
proper morphisms to ${\text {Hot}} ({\bold M})$ which maps $X$ to $W_j
(X) = T(j_X)$ and $f : X \rightarrow Y$ to $W(f)$. If $\bold {j}' =
(j'_X)$ is another choice of compactifications, the maps $W(1_X) :
T(j_X) \rightarrow T(j'_X)$ define an isomorphism of functors from
$W_{\bold {j}}$ to $W_{\bold {j}'}$. In that sense, $X \mapsto W(X)$
is independent of choices and contraviant for proper morphisms.

\medskip

To check that $W(X)$ is covariant for open immersions, notice that
given a compactification $\overline X$ of $X$ there is a contravariant
equivalence of categories between open subschemes of $X$ and closed
subschemes of $\overline X$ containing $\overline X - X$, mapping $(U
\hookrightarrow X)$ to $\left( \overline X - U \hookrightarrow
\overline X \right)$. Since $W(U) = T \left( \overline X -U
\hookrightarrow \overline X \right)$ and $T$ is contravariant on
${\text {Ar}} ({\bold P})$, we conclude that any open immersion $i:U
\hookrightarrow X$ induces $i_* : W(U) \rightarrow W(X)$ and that,
given $i: U \hookrightarrow V$ and $j : V \hookrightarrow X$ two open
immersions, the identity $(ji)_* = j_* \ i_*$ holds. This concludes
the proof of Theorem 2 ii).

\bigskip \subheading{2.3} To prove Theorem 2 iii), let $\bar X$ be a
compactification of $X$, and write $Y=\bar{X}-X$ and $Z=\bar{X}-U$, so
that $T=Z-Y$. Choose non-singular hyper-envelope $\tilde{X}\to
\bar{X}$, $\tilde{Z}\to Z$, and $\tilde{Y}\to Y$ such that there are
maps $$\tilde{Y}. \overset{f}\to \to \tilde{Z}. \overset{g}\to \to
\tilde{X}.$$ lifting the inclusions $Y\to Z \to \bar{X}$. These induce
maps of complexes of motives: $$M(\tilde{X}.) \overset{g^*}\to \to
M(\tilde{Z}.) \overset{f^*}\to \to M(\tilde{Y}.)$$ and a corresponding
triangle of mapping cones: $$ C(g^*)^.[-1]\to C((g\circ f)^*)^.[-1]
\to C(f^*)^.[-1]\to C(g^*)^. $$ which by definition may be rewritten:
$$ W(U)\to W(X) \to W(T) \to W(U)[1]\quad .$$

Property iii) indicates that $W(X)$ behaves like cohomology with
compact supports; see also Theorem 3 below.

\bigskip \subheading{2.4} To prove Theorem 2 i), let $U\subset X$ be a
smooth dense open subset and $T=X-U$ its closed complement. From 2.3
we know that there is a triangle

$$W(T)[-1] \rightarrow W(U) \rightarrow W(X)$$

\noindent in $Hot ({\bold M})$. Assume we know that Theorem 2 i) is
true for $U$ and $T$, i.e. there exist homotopy equivalences $A.\to
W(T)$ and $B.\to W(U)$ where $A_i$ and $B_i$ are motives of varieties
lying in ${\bold V}$ such that ${\dim} (A_i) \leq {\dim} (T)-i$ and
${\dim} (B_i) \leq {\dim} (X)-i$.  In the triangulated category $Hot
({\bold M})$, $W(X)$ is then isomorphic to the cone $C.$ of a map

$$A. [-1] \rightarrow B.$$ (again by \cite{V} II Cor. 1.  ).  Since
$C_i =A_i\oplus B_i$, we have ${\dim} (C_i) \leq {\dim} (X)-i$, and
hence i) is true for $X$.

By noetherian induction we can therefore assume that $X$ is smooth and
quasi-projective.  In that case, let $\overline X$ be a smooth
compactification of $X$ lying in ${\bold V}$, with complement
$T=\overline X -X=D^{\text {red}}$, where $D$ is a divisor with normal
crossing in $\overline X$. Since $W\left(
\overline X
\right)$ is represented by $M\left( \overline X \right)$ in degree
zero, Theorem 2 i) holds for $\overline X$. Using 2.3 again, $W(X)$ is
the cone of a map $$W(\overline X )[-1] \rightarrow W(T)[-1]$$ and,
since $\dim (T)
\leq \dim (X) - 1$, we get i) for $X$ by induction on dimensions.
See 2.7 for a more explicit description of $W(X)$.

\bigskip

\subheading{2.5}
 Let us prove Theorem 2 iv) when $X$ is complete (the general case is
left to the reader). First assume that $X = A \cup B$, where $A$ and
$B$ are closed in $X$. From our discussion in 2.2, diagram (2.2), we
know that there exists a commutative square $$ \CD
\widetilde{A\cap B}.@>v>> \widetilde{B}.\\
    @VuVV @VV{\beta}V \\ \widetilde{A}. @>>\alpha> \widetilde{X}.
\endCD $$
of simplicial varieties in ${\bold V}$ mapping by hyperenvelopes to
the commutative square $$ \CD {A\cap B}@>>> {B}\\ @VVV @VVV \\ {A}
@>>>{X}\ . \endCD $$ The triangle $$ W(X) \rightarrow W(A) \oplus W(B)
\rightarrow W(A\cap B) \rightarrow W(X)[1] $$ is a consequence of the
fact that the total complex of the bicomplex in ${\bold M}$ $$ M
\left(
\widetilde{X}.\right) \overset {\alpha^* \oplus \beta^*}\to
\rightarrow M \left( \widetilde{A}. \right)
 \oplus M\left( \widetilde{B}. \right) \overset {u^* -v^*}\to
\rightarrow M \left( \widetilde{A\cap B}. \right) $$ is
contractible. But this complex is equal, up to a shift, to the mapping
cone of the map of complexes $$ \varphi : S \left( \widetilde{A}.
\hookrightarrow \widetilde{X}. \right) \rightarrow S \left(
\widetilde{A \cap B}. \hookrightarrow \widetilde{B}. \right) .  $$
Since $X-A = B-(A\cap B)$, $\varphi$ is the canonical homotopy between
two representatives of $W(X-A)$ (see 2.3), and this proves our claim.

\medskip

Assume now that the complete variety $X$ is the union of two open
subvarieties $U$ and $V$. Let $A=X-U$ and $B=X-V$. Choose
hyperenvelopes $\widetilde{A}. \rightarrow \widetilde{X}.$ and
$\widetilde{B}. \rightarrow \widetilde{X}.$ of the inclusions $A
\hookrightarrow X$ and $B \hookrightarrow X$. Since $B = U-(U \cap V)$
we know from 2.4 that $W(U \cap V)$ is represented by the complex of
motives $$ \eqalign{ C_1 = & \ C \biggl( C \biggl( M \left(
\widetilde{X}.  \right) \rightarrow M \left( \widetilde{A}.\right)
\biggl)^. [-1] \rightarrow M
\left( \widetilde{B}.\right) \biggl)^. [-1]
\cr = & \ C \biggl( M \left( \widetilde{X}.\right) \rightarrow M
\left( \widetilde{A}.  \right) \oplus M \left( \widetilde{B}.\right)
\biggl)^. [-1]. \cr } $$
 On the other hand, $W(U) \oplus W(V)$ is represented by $$ C_2 = C
\biggl( M \left( \widetilde{X}.\right)
\oplus M \left( \widetilde{X}.\right) \rightarrow M \left(
\widetilde{A}.  \right) \oplus M
\left( \widetilde{B}.\right) \biggl)^.
[-1] , $$ and $W(X)$ is represented by $M \left(
\widetilde{X}.\right)$. The diagonal $M \left( \widetilde{X}.\right)
\rightarrow M \left( \widetilde{X}.\right) \oplus M \left(
\widetilde{X}.\right)$ and the difference map $M \left(
\widetilde{X}.\right) \oplus M \left( \widetilde{X}.\right)
\rightarrow M \left( \widetilde{X}.\right)$ lead to a complex $$ C_1
\rightarrow C_2 \rightarrow M \left( \widetilde{X}.\right) $$ which is
contractible, and this proves the existence of a triangle $$ W(X)[-1]
\rightarrow W(U\cap V) \rightarrow W(U) \oplus W(V) \rightarrow W(X).
$$

  \bigskip \subheading{2.6} The equality $W(X\times Y)=W(X) \otimes
W(Y)$ in $Hot (\bold M)$ follows from the fact that the product of two
envelopes is an envelope. Indeed, we first assume that $X$ and $Y$ are
projective and we let $\widetilde{X}. \rightarrow X$ and
$\widetilde{Y}.  \rightarrow Y$ be two hyperenvelopes. Since coskeleta
commute with products, we deduce from this fact that the product
simplicial scheme $(\widetilde{X} \times
\widetilde{Y}).=(n\mapsto
\widetilde{X}_n \times
\widetilde{Y}_n)$ is an envelope of $X\times Y$. By the
Eilenberg-Zilber theorem (\cite{D-P} 2.9 and 2.16) the associated
complex of motives $M((\widetilde{X}
\times \widetilde{Y}).)$ is homotopy equivalent to the total complex
of the tensor product $M\left( \widetilde{X}.\right) \otimes M\left(
\widetilde{Y}.\right)$, and the equality $W(X\times Y)=W(X) \otimes
W(Y)$ follows.

If we do not assume that $X$ and $Y$ are projective, we let $\overline
X$ and $\overline Y$ be compactifications, $\widetilde{X}.$,
$\widetilde{Y}.$, $\widetilde{S}.$ and $\widetilde{T}.$ be
hyperenvelopes of $\overline X$, $\overline Y$, $S= \overline X -X$
and $T=\overline Y -Y$ respectively, with maps $\widetilde{S}.
\rightarrow \widetilde{X}.$ and $\widetilde{T}. \rightarrow
\widetilde{Y}.$ the cones of which represent $W(X)$ and $W(Y)$. Then
$\overline X \times \overline Y$ is a compactification of $X\times Y$,
and $$\left( \overline X \times \overline Y \right) -(X\times
Y)=\left( S\times
\overline Y \right) \cup \left( \overline X \times T
\right).$$
  By 2.2, diagram (2.2), there is a commutative square of simplicial
varieties in ${\bold V}$ $$ \CD \widetilde{ S \times T.}@>>>
\widetilde{\overline{X} \times T.}\\
    @VVV @VVV \\ \widetilde{S \times \overline{Y}.} @>>>
\widetilde{\overline X \times \overline{Y}.} \endCD $$
 mapping by hyperenvelopes to the square $$ \CD { S \times T}@>>>
{\overline{X} \times T}\\ @VVV @VVV \\ {S \times \overline{Y}} @>>>
{\overline X \times
\overline{Y}} .  \endCD $$
 According to the Mayer-Vietoris property iv), we can represent
$W\left( \left( S \times \overline Y \right) \cup \left( \overline X
\times T \right) \right)$ by the cone of the map $$ M \left(
\widetilde{S \times \overline{Y}.}\right) \oplus M \left(
\widetilde{\overline{X} \times T.}\right) \rightarrow M \left(
\widetilde{ S \times T.} \right), $$ therefore $W(X \times Y)$ is
represented by the total complex of the bicomplex $$ M \left(
\widetilde{\overline X \times \overline{Y}.} \right) \rightarrow M
\left( \widetilde{S \times \overline{Y}.}\right) \oplus M \left(
\widetilde{\overline X \times {T}.} \right) \rightarrow M \left(
\widetilde{S \times T.} \right) .  $$ By the previous step and the
uniqueness proved in 2.3, this bicomplex of motives is homotopy
equivalent to

$$M\left( \widetilde{X}. \right) \otimes M\left( \widetilde{Y}.
\right) \rightarrow M\left( \widetilde{S} \right) \otimes M\left(
\widetilde{Y}.  \right) \oplus M\left( \widetilde{X}. \right) \otimes
M\left( \widetilde{T}. \right) \rightarrow M\left( \widetilde{S}.
\right) \otimes M\left( \widetilde{T} \right),$$

\noindent which is the tensor product of the cone of
 $M\left( \widetilde{X}. \right) \rightarrow M\left( \widetilde{S}.
\right)$ with the cone of $M\left( \widetilde{Y}. \right) \rightarrow
M\left( \widetilde{T}.\right)$.  \qed \bigskip

\bigskip

\subheading{2.7} For future use, we shall give a more precise
description of $W(X)$ when $X$ is smooth and equipped with a
compactification $\overline X$ lying in ${\bold V}$. We assume that
$T=\overline X -X=D^{\text {red}}$, where $D$ is a divisor with normal
crossing in $\overline X$.  Let $Y_1 ,\ldots ,Y_n$ be the irreducible
components of $T$. For any subset $I\subset
\{ 1,\ldots ,n\}$ we let $Y_I = \cap_{i\in I} \ Y_i$ and,
 for any integer $r\geq 1$, we define
$$Y^{(r)} = \coprod_{{\text{card}} (I)=r} \ Y_I \ .$$

\noindent Also we let $Y^{(0)} =Y_{\emptyset} =\overline X$. Clearly
${\dim} \ Y^{(r)} ={\dim} (X)-r$.

\medskip

\noindent If $1\leq k\leq r$ we let $\delta_k :Y^{(r)} \rightarrow
Y^{(r-1)}$ be the disjoint union of the inclusions $Y_I \rightarrow$
\break $Y_J$ where $I$ is the ordered set $\{ i_1 ,\ldots ,i_r \}$ and
$J=\left\{ i_1 ,\ldots ,\widehat{i_k} ,\ldots ,i_r \right\}$.  In
$\Bbb Z \bold V$ we define $$\partial = \sum_{k=1}^{r} (-1)^k \
\delta_k : Y^{(r)} \rightarrow Y^{(r-1)}.$$

\noindent One checks that $\partial \circ \partial =0$ (notice that
$\partial : Y^{(1)} \rightarrow Y^{(0)} =\overline X$ factors via
$T$).

\bigskip

\proclaim{ Proposition 3} The complex

$$M \left( \overline X \right) \ \overset \partial^*\to \to \
M(Y^{(1)}) \ \overset \partial^*\to \to \ M(Y^{(2)}) \ \overset
\partial^*\to \to
 \cdots \ \overset \partial^*\to \to \ M(Y^{({\dim} X)})$$

\noindent is a representative of $W(X)$.

\endproclaim \demo{Proof} Let $\widetilde{Y}. = cosk_0^T
\left( Y^{(1)} \right)$ be the coskeleton of $Y^{(1)}$ over $T$. For
any $r\geq 0$, $\widetilde{Y}_r$ is the disjoint union of the
varieties $Y_{\sigma}$ in ${\bold V}$, where $\sigma$ runs over all
maps $\sigma : \{ 1,\ldots ,r-1 \} \rightarrow \{ 1, \ldots ,n\}$ and
$Y_{\sigma} = Y_{{\text {Im}} (\sigma)}$.  Since the canonical map
$\widetilde{Y}. \rightarrow T$ is a non-singular hyperenvelope,
$M\left( \widetilde{Y}. \right)$ is a representative of $W(T)$ and
$W(X)$ is represented by $C\biggl( M \left( \overline X
\right)
\overset {f^*}\to \rightarrow M \left( \widetilde{Y}. \right) \biggl).
[-1]$, where $f:
\widetilde{Y}. \rightarrow \overline X$ is defined using the
inclusions $Y_{\sigma} \subset \overline X$. Denote by $M.$ the
complex of motives $$ M \left( Y^{(1)} \right) \overset
{\partial^*}\to \rightarrow M \left( Y^{(2)} \right) \overset
{\partial^*} \to \rightarrow \cdots \overset {\partial^*}\to
\rightarrow M \left( Y^{(\dim X)} \right) .  $$ We shall prove
Proposition 3 by exhibiting a homotopy equivalence $\varphi$ from $M
\left( \widetilde{Y}. \right)$ to $M.$. Given any subset $I \subset \{
1,\ldots ,n\}$ of cardinality $r\geq 1$, let $\sigma_I : \{ 1,\ldots
,r\} \rightarrow \{ 1,\ldots ,n\}$ be the unique map of ordered sets
with image $I$. The disjoint union of the identity maps $Y_I
\rightarrow Y_{\sigma_I}$, ${\text {card}} (I) = r$, defines an
inclusion $Y^{(r)} \rightarrow \widetilde{Y}_{r+1}$ and a morphism of
complexes $\varphi : M \left( \widetilde{Y}. \right) \rightarrow M.$.
Now $M \left(
\widetilde{Y}. \right)$ is a representative of $W(T)$ and
$M \left( Y^{(r)} \right)$ represents $W \left( Y^{(r)} \right)$.
  The Mayer-Vietoris property iv)
in Theorem 2 and induction on the number of components of $T$ prove
that $\varphi$ is a homotopy equivalence.  \qed \enddemo \bigskip

\heading 3. Some motivic invariants of varieties \endheading

\medskip

Let $k$ be a field of characteristic zero.  In this paragraph, given a
variety $X$ over $k$, we shall describe several invariants of $X$
which depend only of the associated weight complex $W(X)$.

\bigskip

\subheading{3.1 Weights}

\medskip

\subheading{ 3.1.1} Let $\sim$ be any adequate equivalence relation on
algebraic cycles, ${\bold M}_{\sim}$ be the associated category of
motives (see 1.3), and $\Gamma:{\bold M}_{\sim} \to {\bold A}$ a
covariant (resp.  contravariant) functor from ${\bold M}_{\sim}$ to an
abelian category ${\bold A}$.  If $X$ is any variety over $k$, we may
consider the image $W(X)_{\sim}$ of $W(X)$ in ${\bold M}_{\sim}$ and
the associated complex $\Gamma(W(X)_{\sim})$ in ${\bold A}$. For any
integer $i\geq 0$, we define

$$R^i \Gamma (X) \in Ob ({\bold A}) \quad (\hbox{resp.} \ L_i \Gamma
(X) \in Ob ({\bold A}))$$

\noindent to be the $i$-th cohomology (resp. homology) of $\Gamma
(W(X)_{\sim})$.

\medskip

 From Theorem 2 we conclude that $R^i \Gamma (X)$ is well defined,
contravariant in $X$, and equal to zero if $i> {\dim} (X)$.
Furthermore, when $X$ lies in ${\bold V}$ $R^0 \Gamma (X)=\Gamma
(M(X))$ and $R^i \Gamma (X)=0$ if $i>0$. When $T\subset X$ is a closed
subvariety with complement $U=X-T$, there is a long exact sequence

$$\cdots \rightarrow R^i \Gamma (U) \rightarrow R^i \Gamma (X)
\rightarrow R^i \Gamma (T) \rightarrow R^{i+1} \Gamma (U) \rightarrow
\cdots$$

\noindent Similar properties hold for $L_i F$.

\bigskip

\subheading{ 3.1.2} A basic example comes when $k={\Bbb C}$, $\sim$
is homological equivalence, and $\Gamma (X,p) =p_* H^n (X({\Bbb C}),
A)$ is the singular (= Betti) cohomology of the motive $M=(X,p)$ with
constant coefficients in a given ring $A$, for a fixed value of the
integer $n\geq 0$, correspondences acting on cohomology in the usual
way.  Clearly $\Gamma$ defines a contravariant functor $H^n$ from
${\bold M}_{\sim}$ to the category of finitely generated $A$-modules.

\medskip
 Choose a compactification $\overline X$ of $X$, let $j:Y=\overline X
- X\rightarrow \overline X$ be its complement and let $\widetilde
j:\widetilde Y.
\rightarrow \widetilde X.$ be
 a non-singular hyperenvelope of $j$. By definition $W(X)$ is
represented by $S(\widetilde j)$. Therefore

$$R^i H^n (X)=H^i \left( * \mapsto H^n (S(\widetilde j)^*
,A)\right).$$

\noindent Since envelopes are proper and surjective they define
morphisms ``de descente cohomologi\-que universelle'' by \cite{D} ,
5.3.5 (II). Therefore, the hypercohomology of $S(\widetilde j)$ is the
relative cohomology $H^* \left( \overline X ,\overline Y ,A\right)
=H_c^* (X({\Bbb C}),A)$, and $R^i H^n (X)$ coincides with the term
$E_2^{i,n}$ of the corresponding descent spectral sequence, which is
thus independent of choices.

When $A={\Bbb Q}$, this weight spectral sequence degenerates at $E_2$
(op.cit. Proposition (8.1.20)) and we get (loc.cit.)

$$H^i \left( * \mapsto H^n (S(\widetilde j)^*,{\Bbb Q})\right) =
E_2^{i,n} =E_{\infty}^{i,n} ={\text gr}^W_n \ H_c^{i+n} (X({\Bbb C}),
{\Bbb Q}).$$ In other words, using Theorem 2, we have the following
result
\bigskip

\proclaim { Theorem 3}
 The cohomological descent spectral sequence $$E_2^{i,n}\Rightarrow
H_c^{i+n} (X({\Bbb C}),A)$$ is independent of choices when it comes
from hyperenvelopes as above. It defines a canonical increasing weight
filtration $F_n^W H_c^k (X({\Bbb C}),A)$ on the cohomology with
compact support of $X({\Bbb C})$ with constant coefficients in the
ring $A$.  This filtration has length at most $dim(X)+1$.  It is
compatible with products, pull-back by proper maps, and push-forward
by open immersions.

When $A={\Bbb Q}$, the filtration coincides with the weight filtration
defined by Deligne in \cite{D} and $$R^i H^n(X)\otimes {\Bbb Q}
={\text gr}_n^W H_c^{i+n} (X({\Bbb C}),{\Bbb Q}),$$ where ${\text
gr}_n^W$ is the subquotient of weight $n$.

\endproclaim \bigskip

When $X$ is projective and $k > 0$, $F_{k-1}^W H_c^{k} (X({\Bbb
C}),A)$ is the kernel of the map $$\pi^*:H_c^{k} (X({\Bbb C}),A)\to
H_c^{k} (X'({\Bbb C}),A)$$ for any resolution of singularities $\pi :
X' \to X$.  That this kernel is independent of choices was first
observed by Grothendieck \cite{Gr2}.

\medskip

Theorem 5.13 in \cite{G-N} leads similarly to a canonical weight
filtration on the cohomology
without support and with
 arbitrary coefficients.

\bigskip

\subheading{ 3.1.3}
In order to check that the weight filtration on $H_c^{k} (X({\Bbb
C}),{\Bbb Z})$ is non-trivial and cannot be defined in a simple way
from the cohomology with rational coefficients, let us consider the
following example.

 Let $T$ be an abelian surface, and let $i:T\to T$ be the involution
defined by $i(x)=-x$. The quotient surface $S=T/<1,i>$ is projective
with sixteen ordinary double points $\{p_1,...,p_{16}\}$ of type
$A_1$.  On resolving these singularities, we obtain the Kummer surface
$\tilde{S}$ associated to $T$, which is a $K_3$-surface. If $\pi:
\tilde{S}\to S$ is the resolution, the sixteen exceptional curves
$W_i=\pi^{-1}(p_i)$ are all isomorphic to ${\Bbb P}^1$. We may view
$S$ as the pushout in the following diagram $$\CD
\coprod E_i@>>>\tilde{S}\\ @VVV @VVV \\ \coprod p_i @>>> S.
\endCD $$
In this case one may check that the weight spectral sequence of
Theorem 3 is simply the corresponding Mayer-Vietoris exact sequence:
$$\ldots\to H^n(S,{\Bbb Z})\to H^n(\tilde{S},{\Bbb Z})\oplus H^n(
\coprod p_i,{\Bbb Z})
\to H^n(\coprod E_i,{\Bbb Z})\to  H^{n+1}(S,{\Bbb Z})\to \ldots$$
In particular we find that the sequence $$0\to H^2(S,{\Bbb Z})\to
H^2(\tilde{S},{\Bbb Z})
\to {\Bbb Z}^W\to  H^{3}(S,{\Bbb Z})\to 0$$
is exact, where $W=\{ W_1,...,W_6 \}$, so that $$H^3(S,{\Bbb
Z})=F_2^W\ H^3(S,{\Bbb Z})= Gr_2^W\ H^3(S,{\Bbb Z})$$ while
$H^1(S,{\Bbb Z})=0$.

In Proposition 5.5 of Chapter VIII of \cite{B-P-V} , it is shown that
the sublattice of $H_2(S,{\Bbb Z})$ generated by the classes of the
divisors $W_i$ is of index $32$ in the smallest primitive sublattice
containing it, with the quotient being 2-torsion; dualizing we find $$
Gr_2^W\ H^3(S,{\Bbb Z})= H^3(S,{\Bbb Z}) = ({\Bbb Z}/2{\Bbb Z})^5.$$
On the other hand, if $X$ is an Enriques surface (op. cit.  Chapter
VIII, Section 15), $\pi_1(X)={\Bbb Z}/2{\Bbb Z}$, and so by Poincar\'e
duality $ H^3(X,{\Bbb Z}) = {\Bbb Z}/2{\Bbb Z}$ and $H^1(X,{\Bbb Z}) =
0$. Since $X$ is smooth and projective $$Gr_3^W\ H^3(X,{\Bbb Z})=
H^3(X,{\Bbb Z}) = {\Bbb Z}/2{\Bbb Z}.$$ Taking the cartesian product
$S\times X$ we obtain a four-fold with $$ H^3(S\times X,{\Bbb Z}) =
H^3(S,{\Bbb Z})\oplus H^3(X,{\Bbb Z}) = ({\Bbb Z}/2{\Bbb Z})^6$$
(since $H^1(S,{\Bbb Z}) = H^1(X,{\Bbb Z}) = 0$ and $H^2(S,{\Bbb Z})$
is torsion free), and $$Gr_3^W\ H^3(S\times X,{\Bbb Z})= {\Bbb
Z}/2{\Bbb Z}$$ $$Gr_2^W\ H^3(S\times X,{\Bbb Z})= ({\Bbb Z}/2{\Bbb
Z})^5.$$

\bigskip

\subheading{ 3.1.4} When $M=(X,\pi)$ is a Chow motive and $p\geq 0$ an
integer, we may consider the Chow cohomology (resp. homology) group
$\pi_* C H^p (X)$ (resp. $\pi_* C H_p (X)$ ).

\bigskip

\proclaim{ Proposition 4} When $X$ is a (possibly singular) complete
variety over $k$, the group $L_0 C H_p (X)$ coincides with $C H_p (X)$
and $R^0 CH^p (X)$ is the operational Chow group $A^p( X \ \overset
{\text id}\to \to X )$ (see \cite{F} ).
\endproclaim \bigskip

\demo{ Proof} This follows from the descent properties of Chow
homology and from a result of Kimura \cite{Ki}, see \cite{B-G-S}
Appendix.  \qed \enddemo \medskip

More generally, if $Z.$ is a non-singular hyperenvelope of a
projective variety $X$, for any $p\geq 0$ there is a canonical weight
spectral sequence

$$E_2^{st}(X) =H^s \left( n \mapsto H^t (Z_n ,{\Cal K}_p )\right)
\Rightarrow H^{s+t} (Z. ,{\Cal K}_p)$$ converging to the
hypercohomology of the simplicial scheme $Z.$ with coefficients in the
Zariski sheaf ${\Cal K}_p$.  Up to canonical isomorphism, this
spectral sequence is independent of the choice of $Z.$ from $E_2$ on
(indeed, a map of hyperenvelopes induces a morphism of spectral
sequences and $E_2^{st}(X) $ depends only on $W(X)$, so we can argue
as in 2.2).  From Proposition 3 and the Gersten conjecture \cite{Q} we
conclude that $$R^0 CH^p (X)=E_2^{0,p}(X) .$$

\bigskip

\subheading{ 3.1.5}
 Assume $X$ is a smooth variety and that $\overline X$ is a
compactification of $X$ lying in ${\bold V}$, with complement
$T=D^{\text{red}}$, where $D$ is a divisor with normal crossing in
$\overline X$ as in 2.7.  From Proposition 3 it follows that $R^i CH^p
(X)$ is the $i$-th cohomology of the complex

$$CH^p \left( \overline X \right) \ \overset {\partial^*}\to \to \
CH^p \left( Y^{(1)} \right) \ \overset {\partial^*}\to \to \ CH^p
\left( Y^{(2)} \right)\ \overset {\partial^*}\to \to \cdots\leqno
(3.1)$$

\noindent when $L_i CH_p (X)$ is the $i$-th homology of

$$\cdots \overset {\partial^*}\to \to \ CH_p \left( Y^{(2)} \right) \
\overset {\partial_*}\to \to \ CH_p \left( Y^{(1)} \right) \ \overset
{\partial_*}\to \to \ CH_p \left( \overline X \right),$$

\noindent where $\partial^* =\Sigma (-1)^k \ \delta_k^*$ and
$\partial_* =\Sigma (-1)^k \ \delta_{k_*}$ (with the same notation as
in 2.7).  So these groups are independent of the choice of the
compactification $\overline X$..

\medskip

Let $i:Y^{(1)}\rightarrow \overline X$ be the disjoint union of the
inclusions $Y_i\to \overline X, 1\leq i\leq n $, $CH^p (D)=A^p \left(
D \ \overset {\text id}\to \to \ D \right)$ the operational Chow group
of $D$, $d= {\dim} (X)$ and

$$i^* i_* :CH_{d-p} \left( Y^{(1)} \right) \rightarrow CH^p \left(
Y^{(1)} \right)$$

\noindent the composite of the maps

$$CH_{d-p} \left( Y^{(1)} \right) \ \overset {i_*} \to \to \ CH_{d-p}
\left( \overline X \right)=CH^p \left( \overline X \right) \ \overset
{i^*} \to \to CH^p \left( Y^{(1)} \right) .$$

\noindent It follows from the result above that the cohomology of the
${\Bbb Z}$-graded complex

$$ \cdots \rightarrow CH_{d-p} \left( Y^{(2)} \right) \ \overset
{\partial_*}\to \to \ CH_{d-p} \left( Y^{(1)} \right) \ \overset {i^*
i_*}\to \to \ CH^p
\left( Y^{(1)} \right) \rightarrow
CH^p \left( Y^{(2)} \right) \rightarrow \cdots$$

\noindent is independent of choices. In particular we recover Theorem
2.2.1 of \cite{B-G-S}. If we follow the analogy of that paper with
K\"ahler geometry, we may view the complex above as analogous to the
following one, defined for any smooth projective complex manifold $M$:

$$ \cdots \ \overset {d}\to \to \ A^{p-3,p-1} (M) \oplus A^{p-2,p-2}
(M) \oplus A^{p-1,p-3} (M) \ \overset {d} \to \to \ A^{p-2,p-1} (M)
\oplus A^{p-1,p-2} (M)$$ $$\overset {d}\to \to \ A^{p-1,p-1} (M)
\overset {dd^c}\to \to \ A^{p,p} (M) \ \overset {d} \to \to \
A^{p+1,p} (M) \oplus A^{p,p+1} (M) \ \overset {d} \to \to\ \cdots , $$
\noindent where $A^{pq} (M)$ are complex forms of type $(p,q)$. This
complex is known to compute the Deligne cohomology groups $H_D^*
(M,{\bold R} (p))$ when appropriate reality conditions are imposed on
forms (see \cite{B} Theorem 1.10). We conclude from this that $R^i
CH^p (X)$ is analogous to $H_D^{2p+i-1} (M,{\bold R} (p))$ and $L_i
CH_{d-p} (X)$ is analogous to $H_D^{2p-i} (M,{\bold R} (p))$ when
$i\geq 1$.

\bigskip

\subheading{ 3.1.6} Assume now that $\sim$ is numerical equivalence.
It was shown by Jannsen \cite{J} that ${\bold M}_{\text{num}}$ is an
abelian semi-simple category. If we apply the discussion of 3.1.1 to
the identity functor, we get canonical motives $W^i (X)$ attached to
$X$, namely the cohomology groups of $W(X)_{\text{num}}$. These
satisfy the properties of 3.1.1 as well as the K\"unneth formula

$$W^i (X\times Y) = \ \oplus_{j+k=i}^{} \ W^j (X) \otimes W^k (X).$$

\noindent If one knew that numerical equivalence implies homological
equivalence, for any $n\geq 0$, the group ${\text gr}^W_i H_c^{n+i}
(X({\Bbb C}), {\Bbb Q})$ would depend only on $W^i(X)$ (when the
ground field is ${\Bbb C}$).

\bigskip

\subheading{ 3.2. The Grothendieck group of motives}

\medskip

\subheading{ 3.2.1}
 Let $\sim$ be any adequate equivalence relation on algebraic cycles
and ${\bold M}_{\sim}$ the associated category of motives .  The
Grothendieck group of this category is the quotient $K_0({\bold
M}_{\sim})$ of the free abelian group on the isomorphism classes $[M]$
of objects $M$ in $\bold M$ by the subgroup generated by elements of
the form $[M]-[M']-[M'']$ whenever $M\simeq M'\oplus M''$.  On the
other hand, if $Hot^b({\bold M}_{\sim})$ denotes the category of
bounded cochain complexes in ${\bold M}_{\sim}$ up to homotopy, we may
consider its Grothendieck group $K_0(Hot^b({\bold M}_{\sim}))$, which
is generated by objects in $Hot^b({\bold M}_{\sim})$ with the relation
$[Y^.]=[X^.]+[Z^.]$ whenever there exists a triangle $$X^.\to Y^. \to
Z^.
\to X^.[1].$$ \bigskip

\proclaim{ Lemma 3} The obvious functor ${\bold M}_{\sim}\to
Hot^b({\bold M}_{\sim})$ induces a group isomorphism $K_0({\bold
M}_{\sim})\to K_0(Hot^b({\bold M}_{\sim}))$.

\endproclaim \medskip

\demo{Proof} This fact is true for any pseudo-abelian category instead
of ${\bold M}_{\sim}$.  Indeed, when $M\simeq M'\oplus M''$ there is a
triangle $$M'\to M \to M'' \to M'[1],$$ \cite{V} Cor. 1.2.3, so we get
a morphism $\phi:K_0({\bold M}_{\sim})\to K_0(Hot^b({\bold
M}_{\sim}))$.

Given a bounded complex $M.$ in ${\bold M}_{\sim}$ we let $$\chi (M^.)
= \sum_i (-1)^i [M^i] \in K_0({\bold M}_{\sim}).$$ When $f:M^.\to N^.$
is a morphism of complexes and $C(f)^.$ is its mapping cone, we have
$$\chi (C(f)^.) = \chi (N^.)- \chi (M^.).$$ So to prove that $\chi$
induces a morphism from $K_0(Hot^b({\bold M}_{\sim}))$ to $K_0({\bold
M}_{\sim}) $ all we need to check is that $\chi (M^.) =0$ when $M^.$
is contractible.  For this we proceed by induction on the length $k$
of $M^.$. Let $h:M^{i+1}\to M^i$ be such that $$dh+hd = id_{M^.}.$$ It
follows that $(hd)^2=hd$. Let $A^{k-1}$ be the image of the projector
$hd$ in $M^{k-1}$ and $B^{k-1}$ its complement.  The map $d:M^{k-1}\to
M^k$ is zero on $B^{k-1}$ and induces an isomorphism from $A^{k-1}$ to
$M^k$.  Let $N^.$ be the complex obtained from $M^.$ by replacing
$M^{k-1}$ by $B^{k-1}$ and $M^k$ by zero.  We get $$\chi (M^.) = \chi
(N^.) + [A^{k-1}] - [X^k] = \chi (N^.).$$ Since $N^.$ is contractible,
this proves that $\chi (M^.) =0$ by induction.  It is now easy to
check that $\phi$ and $\chi$ are isomorphisms inverse to each other.
\qed\enddemo

Using this lemma and Theorem 2 we get

\bigskip

\proclaim{ Theorem 4} Any quasi-projective variety $X$ has a class
$[X]\in K_0 ({\bold M}_{\sim})$ characterized by the following
properties:

\smallskip

\item{i)} If $X$ lies in ${\bold V}$, $[X]$ is the class of the motive
$(X,1_X)$;

\smallskip

\item{ii)} If $Y\subset X$ is a closed subvariety in $X$,

$$[X]=[Y]+[X-Y].$$ \endproclaim \medskip

\demo{Proof} Define $[X]$ to be the class of $W(X)$ in $K_0
(Hot^b({\bold M}_{\sim}))=K_0 ({\bold M}_{\sim})$.  Then property i)
is clear and ii) follows from Theorem 2 iii).  To see that $[X]$ is
uniquely characterized by i) and ii) we proceed by induction on
dimensions.  If $U\subset X$ is a smooth dense open subset and
$\overline U$ a smooth compactification of $U$ lying in ${\bold V}$
(which exist by resolution of singularities), we get
$$[X]=[U]+[X-U]=[\overline U]-[\overline U-U]+[X-U],$$ which fixes
$[X]$ uniquely.  \qed\enddemo Theorem 4 answers positively a question
of Serre, \cite{Se} p.341.
\bigskip

\subheading{ 3.2.2} \medskip

\proclaim{ Proposition 5} The class $[X]\in K_0 ({\bold M}_{\sim})$
has the following properties:

\smallskip

\item{i)} If $U$ and $V$ are two locally closed
 subvarieties in $X$, then $$[U\cap V]+[U\cup V] = [U] + [V];$$
\smallskip

\item{ii)} If $f:X\rightarrow B$ is a fibration of fiber $F$ which is
locally trivial for the Zariski topology of $B$, then

$$[X]=[F]\cdot [B].$$

\endproclaim \medskip

\demo{Proof} For i) notice that $V-(U\cap V) = (U\cup V)-V$,
 therefore, if $U$ and $V$ are open or if $U$ and $V$ are closed, we
deduce from ii) in Theorem 4 (or from Theorem 2 iv)) that $$[V]-[U\cap
V] = [V-(U\cap V)] = [U\cup V] - [V].$$ The general case follows from
these two.

To prove ii), cover $B$ by a finite collection of open subsets
$U_{\alpha}, \alpha\in A,$ such that $f$ is trivial over $U_{\alpha}$.
Using i) and induction on the cardinality of $A$ one is reduced to
proving ii) when $f$ is trivial, i.e.  $$[X\times Y]=[X] \cdot [Y]
\quad \hbox{in} \quad K_0 ({\bold M}_{\sim}).$$ Such an equality
follows from Theorem 2 v) or can be deduced directly from the case
where $X$ and $Y$ are smooth and projective by induction on $\dim
(X)+\dim (Y)$, as in the proof of Theorem 4 above.  \qed\enddemo
\bigskip

\subheading{ 3.2.3} It follows from Theorem 4
 that the class of the affine line in $K_0 ({\bold M}_{\sim})$ is
class of the Tate motive $\Bbb L$: $$[{\Bbb A}^1] = [{\Bbb P}^1] - [1]
=[{\Bbb L}].$$

If $X$ is the affine cone with base a variety $Y$ lying in ${\Bbb V}$,
then, as in \cite{Se}, $$[X] = [1]+[Y]\otimes [{\Bbb L}] - [Y].$$

Much more elaborate cases can be found in the recent paper by Manin
\cite{M2}.

 \bigskip

\subheading{ 3.2.4}
 Theorem 4 is related to the theory of mixed motives as follows. In
\cite{Vo}, Voevodsky associates to any perfect field $k$ a
triangulated category ${\bold {DM}}_{\text {gm}}^{\text {eff}}$ of
``effective geometrical motives over $k$''.  When ${\text {char}}(k) =
0$, any variety $X$ of finite type over $k$ has classes $M(X)^c$ and
$M(X)$ in ${\bold {DM}}_{\text {gm}}^{\text {eff}}$.  Their properties
are listed in \cite{Vo} 2.2. In particular, $M(X) = M(X)^c$ when $X$
is smooth and proper over $k$, and the restriction of $M$ to ${\bold
V}$ can be factorized through an additive functor from the category of
Chow motives ${\bold M}$ to ${\bold {DM}}_{\text {gm}}^{\text {eff}}$.
Therefore we get a group morphism $$ \varphi : K_0 ({\bold M})
\rightarrow K_0 \left({\bold {DM}}_{\text {gm}}^{\text {eff}}\right) ,
$$ from $K_0 ({\bold M})$ to the Grothendieck group of V\oe vodsky's
triangulated category.

For any $X$ of finite type over $k$, $\varphi ([X])$ is the class of
$M(X)^c$.  Indeed, this is true by definition when $X$ lies in ${\bold
V}$, and the general case follows using Theorem 4 ii) and Property 2
of $M(X)^c$ in \cite{Vo} 2.2.

It seems a hard problem to decide whether $\varphi$ is a group
isomorphism or not.

\bigskip

\subheading{ 3.3.  Numerical invariants}

\medskip

\subheading{ 3.3.1} From Theorem 4 it follows that any additive map $
h :Ob ({\bold M}_{\sim}) \to A$, where $A$ is an abelian group defines
for each variety $X$ over $k$ a class $h(X)\in A$.  This class has the
properth that $h(X)=h(Y)+h(X-Y)$ when $Y$ is closed in $X$ .  \medskip

For instance, if $k={\Bbb C}$, $\sim$ is homological equivalence,
$A={\Bbb Z}$, $n\geq 0$ is a fixed integer and $$h(M)= \dim_{\Bbb Q}\
p_*H_c^{n} (X({\Bbb C}),{\Bbb Q})$$ if $M=(X,p)$, we deduce from
Theorem 3 that, for an arbitrary variety $X$, $h(X)\in {\Bbb Z}$ is
the "$n$-th virtual Betti number" $$h^n(X)=\sum_{i\geq 0} \ (-1)^i
\ {\dim_{\Bbb Q}} \ {\text gr}^W_n H_c^{n+i} (X({\Bbb C}),{\Bbb Q})\in
{\Bbb Z}.$$

A much finer invariant than the virtual Betti number, which includes
torsion data, may be found by taking $A$ to be the Grothendieck group
of all finitely generated abelian groups with respect to {\it direct
sum}.  Thus:
 $$A={\Bbb Z}\oplus \bigoplus_{p }Y,$$ where $p$ runs
over all prime integers and $Y$ is the free abelian group on the set
$\Bbb N$ of natural numbers.  That is, Y is the group of functions
$\phi:{\Bbb N}\to {\Bbb Z}$ such that $\phi(n)=0$ for all but finitely
many $n$.

If $G$ is an abelian group the corresponding class $h(G)\in A$ is $$
{\text rk}(G)\oplus \bigoplus_p \phi_p$$
 where

$$ G\simeq {\Bbb Z}^{{\text rk} G}\oplus \bigoplus_p ({\Bbb
Z}/p^n)^{\phi_p(n)}$$
 Thus from the class $h(G)$ we recover the
{\it isomorphism} class of $G$.  Now we may take, for $X$ a smooth
projective variety, $h(X)$ to be the class of $H^n(X,{\Bbb Z})$ in
this class group.  The corresponding invariant of singular varieties
seems to be new, and combines information about the torsion in
integral cohomology with the weight structure.  From this invariant
one can deduce several numerical invariants. For example the invariant
obtained by taking $\text {card} (\text {tors}(H^n(X,{\Bbb Z})))$ for
$X$ smooth and projective, which was first suggested by Totaro.

Given $p,q\geq 0$ two integers, we may also consider algebraic De Rham
cohomology and, for any motive $M=(X,p)$, the integer $$h(M)=
{\dim_k}p_*H^p(X,\Omega^q),$$ where correspondences act via their De
Rham fundamental class.  The corresponding invariant $h^{p,q}(X)\in
{\Bbb Z}$ (for $X$ any variety over $k$) is mentioned by Grothendieck
in \cite{Gr3} p.191.  In terms of the mixed Hodge structure on
cohomology \cite{D} one gets, if $k\subset {\Bbb C}$ say,
$$h^{p,q}(X)= \sum_{i\geq 0}\ (-1)^i \ {\dim_{\Bbb C}}
\ {\text gr}_F^q\ {\text gr}_{p+q}^W
 H_c^{p+q+i}(X({\Bbb C}),{\Bbb C}).$$ It follows from the axiomatic
definition of these numbers that, for every $X$, one has
$$h^n(X)=\sum_{p+q=n} h^{p,q}(X)$$ and $h^{p,q}(X)=h^{q,p}(X)$.
\bigskip

\subheading{ 3.3.2} If $\sim$ is numerical equivalence and $k={\Bbb
C}$ the usual Euler characteristic

$$\chi (X)=\sum_{n\geq 0} \ (-1)^n \ {\dim} \ H^n (X({\Bbb C}),{\Bbb
Q})$$

\noindent extends to all motives by the formula

$$\chi (M)=\Delta \cdot p$$ where $M=(X,p)$, $\Delta \in End_{\bold C}
(X)$ is the diagonal and $\Delta \cdot p$ denotes its intersection
number with $p$ in the Chow ring of $X\times X$.  From Theorem 3 and
3.1.5 we get

$$\chi (W^i X)=\sum_{n\geq 0} \ (-1)^n \ {\dim} \ {\text gr}^W_i
H_c^{n+i} (X({\Bbb C}),{\Bbb Q})$$

\noindent for any complex variety $X$.

\bigskip

\heading 4. Blow ups \endheading

\medskip

\subheading{ 4.1} For any integer $m\geq 0$ and any noetherian
(separated) scheme $S$, denote by $K_m (S)$ the higher $K$-group of
perfect complexes on $S$ ( \cite{T-T} 3.1) and by $K'_m (S)$ the
higher $K$-group of coherent $\Cal O_S$-modules on $S$ \cite{Q}.  If
$S$ is noetherian and regular and if $X/S$ is any scheme of finite
over $S$ we denote by $A_p (X/S)$ the homological Chow group of cycles
of relative dimension $p$ on $X$ and by $A^p (X/S)$ the cohomological
(operational) Chow group, as in \cite{F}, Appendix.

\medskip

Now let $X$ be a noetherian scheme, $i:Y\rightarrow X$ a closed
immersion.  Consider the blow up $f:X' \rightarrow X$ of $X$ along
$Y$, and let $Y'$ be the inverse image of $Y$ so that we get a
cartesian square of proper maps

 $$ \CD Y'@>j>> X'\\ @VgVV @VVfV \\ Y @>i>> X\ .  \endCD $$

\medskip

\proclaim { Theorem 5}

\smallskip

\item{i)} Assume $i$ is regular. Then, for any $m\geq 0$, the
following sequence is exact: $$0\rightarrow K_m (X) \ \overset (i^*
,-f^*)\to \to \ K_m (Y) \oplus K_m (X') \ \overset {g^* +j^*} \to \to
\ K_m (Y') \rightarrow 0. \leqno (4.1)$$

\item{ii)} Assume $i$ is regular and there exists an ample line bundle
on $X$.  Then, for any $m\geq 0$, the following sequence is exact
$$0\rightarrow K'_m (Y') \ \overset (g_* ,j_*)\to \to \ K'_m (Y)
\oplus K'_m (X') \ \overset{i_* -f_*}\to \to \ K'_m (X) \rightarrow 0.
\leqno (4.2)$$

\item{iii)} Assume that $X$ is of finite type over a regular
noetherian scheme $S$ and that $i$ is regular. Then, for any $p\geq
0$, the following sequences are exact: $$0\rightarrow A_p (Y'/S) \
\overset (g_* ,j_*)\to \to \ A_p (Y/S) \oplus A_p (X'/S) \
\overset{i_* -f_*}\to \to \ A_p (X/S) \rightarrow 0 \leqno (4.3)$$
\noindent and $$0\rightarrow A^p (X/S) \ \overset (i^* ,-f^*)\to \to \
A^p (Y/S) \oplus A^p (X'/S) \ \overset {g^* +j^*}\to \to \ A^p (Y'/S)
\rightarrow 0.  \leqno (4.4)$$

\item{iv)} Assume that $S={\text {Spec}}(k)$
 where $k$ is any field of characteristic zero, and that $X$ is smooth
and projective over $k$. Then, for any closed immersion $i$ (regular
or not), the sequences (4.3) and (4.4) are exact.  \endproclaim
\bigskip

\subheading{ 4.2} In order to prove Theorem 5, we first remark that in
cases i), ii) and iii) the closed immersion $i:Y\rightarrow X$ is
regular, therefore the same is true for $j:Y' \rightarrow X'$, and the
maps $g$ and $f$ are locally complete intersection morphisms.
Consequently all the maps $i$, $j$, $f$ and $g$ are perfect morphisms
and in particular of finite Tor-dimensions. It follows that together
with the morphisms $i^*$, $j^*$, $f^*$, $g^*$ on $K_m$ and $A^p$ there
are push-forward morphisms $i_*$, $j_*$, $f_*$, $g_*$ (which do not
preserve the degree), and similarly $K'_m$ and $A_p$ are both
covariant and contravariant for these maps (for $K'_m$ this uses the
fact that $X$ has an ample line bundle). If $N$ is the normal sheaf of
$Y$ in $X$, $Y'$ is the projective bundle ${\bold P} (N)$ (in the
sense of Grothendieck).  Let $F={ \ker} (g^* N \rightarrow {\Cal O}
(1))$ and $\Lambda^i F$ the exterior powers of this locally free sheaf
on $Y'$.

\medskip

To prove Theorem 5 i) we first consider the map $g^! :K_m (Y)
\rightarrow K_m (Y')$ defined by $g^! (x)=\lambda_{-1} (F) \ g^* (x)$.
Since $f_* f^* (x) =f^* f_* (1) x=x$ for any $x\in K_m (X)$ (
\cite{Gr1} VII Lemme 3.5), the morphism $f^*$ is injective, and since

$$g_* \ g^! (x)=g_* \ (\lambda_{-1} (F)) \ x=x$$

\noindent ( \cite{Gr1}
 VI.5.9 or (6) below), the exactness of (4.1) is equivalent to saying
that $j^*$ induces an isomorphism

$${K_m (X') \over f^* \ K_m (X)} \rightarrow {K_m (Y') \over g^* \ K_m
(Y)} \ .$$

\noindent According to \cite{T2} (2.2.1) the following sequence is
exact

$$0 \rightarrow K_m (Y) \ \overset (g^! ,-i_*)\to \to \ K_m (Y')
\oplus K_m (X) \ \overset{j_* +f^*}\to \to \ K_m (X') \rightarrow 0.$$

\noindent Therefore $j_* : K_m (Y') \rightarrow K_m (X')$ induces an
isomorphism

$${K_m (Y') \over g^! \ K_m (Y)} \rightarrow {K_m (X') \over f^* \ K_m
(X)}$$

\noindent and all we need to show is that $j^* \ j_*$ induces an
isomorphism

$${K_m (Y') \over g^! \ K_m (Y)} \rightarrow {K_m (Y') \over g^* \ K_m
(Y)} \ . \leqno (4.5)$$

\noindent Now ${\Cal O} (1)$ is the conormal bundle of $Y'$ in $X'$,
and the composite $j^* \ j_*$ is the product by the element $\xi =1-L
\in K_0 (Y')$, with $L=[{\Cal O} (-1)]$ being the class of the dual of
${\Cal O} (1)$ ( \cite{T2} (3.1.4)). Any element $x$ in $K_m (Y')$ can
be written uniquely as a sum

$$x=\sum_{i=0}^{d-1} g^* (y_i ) \ L^i$$

\noindent where $y_i \in K_m (Y)$ and $d$ is the rank of $N$, since
$Y' ={\bold P} (N)$. It follows that $x$ can be written uniquely as

$$x=\sum_{i=0}^{d-1} g^* (z_i ) \ \xi^i$$

\noindent with $z_i \in K_m (Y)$. Therefore $x-g^* (z_0)$ is a
multiple of $\xi$, and the map (4.5) is surjective.

\medskip

To prove that (4.5) is injective we first get from the definition of
$F$ the formulae

$$[g^* \Lambda^i N] = [\Lambda^i F] + [\Lambda^{i-1} F] \ [{\Cal O}
(1)] \quad , \quad i\geq 1,$$

\noindent in $K_0 (Y')$. By induction on $i$ it follows that

$$[\Lambda^i F] = (-1)^i \ [{\Cal O} (1)]^i +\sum_{j<i} \alpha_j \
[{\Cal O} (1)]^i$$

\noindent with $\alpha_j \in g^* K_0 (Y)$ and

$$\lambda_{-1} (F)=\xi^{d-1} +\sum_{i=0}^{d-1} \beta_i \ \xi^i \leqno
(4.6)$$

\noindent with $\beta_i \in g^* K_0 (Y)$. Therefore any element $x$ in
$K_0 (Y')$ can be written uniquely as

$$x=\lambda_{-1} (F) \ g^* (y) +\sum_{i=0}^{d-2} g^* (y_i ) \ \xi^i$$

\noindent with $y,y_i \in K_m (Y)$. If we assume that $\xi x =g^* (z)$
we get

$$\sum_{i=0}^{d-2} g^* (y_i ) \ \xi^{i+1} +g^* (y \ \lambda_{-1}
(N)-z)=0$$

\noindent (since $\xi \lambda_{-1} (F)=g^* \ \lambda_{-1} (N)$),
therefore $y_i =0$, $0\leq i\leq d-2$ and $x=\lambda_{-1} (F) \ g^*
(y)$. This proves that (4.5) is injective.

\bigskip

\subheading{4.3} To prove Theorem 5 ii) we first notice that the open
complement $U=X'-Y'$ is isomorphic to $X-Y$. From \cite{Q} \S 7
Proposition 3.2 and Proposition 2.11 we get a commutative diagram $$
\CD
 K'_{m+1} (U) @>>> K'_m (Y') @>j_*>> K'_m (X') @>>> K'_m (U) @>>>
K'_{m-1} (Y')\\ @VVV @VVg_*V @VVf_*V @VVV @VVV\\ K'_{m+1} (U) @>>>
K'_m (Y) @>i_*>> K'_m (X) @>>> K'_m (U) @>>> K'_{m-1} (Y) \endCD $$
\noindent Therefore the sequence (4.2) is exact in the middle and all
we need to show is that $f_* :K'_m (X')
\rightarrow K'_m (X)$ is surjective. This is where we use that $X$ has
an ample line bundle, since we get then a morphism $f^* :K'_m (X)
\rightarrow K'_m (X')$ such that $f_* \ f^* (x)=x \ f_* (1)$ (
\cite{Q} Proposition 2.10). Since $f_* (1)=1$ in $K_0 (X)$ (
 \cite{Gr1} VII Lemme 3.5), this concludes the proof.

\bigskip

\subheading{ 4.4} We prove Theorem 5 iii) by using the fact that the
groups $A_p (X/S)$ and $A^p (X/S)$ satisfy all the properties stated
in \cite{F} when $S={\text {Spec}}(k)$, except those involving
external products (
\cite{F} 20.1).

\medskip

To simplify the notations, we write $A_p (X)$, $A_p (Y)$, $A^p (X)
\ldots$
 instead of $A_p (X/S)$, $A_p (Y/S)$, $A^p (X/S) \ldots$.

\medskip

We first prove that (4.3) is exact. Let $c_{d-1} (F) \in A^{d-1} (Y')$
be the top Chern class of $F$ and $g^! :A_p (Y) \rightarrow A_p (Y')$
be the map sending $x$ to $g^* (x) \ c_{d-1} (F)$. According to
\cite{F} Proposition 6.7 we have $f_* \ f^* ={\text id}$, and, by
 \cite{F} Example 3.3.3, $g_* \ g^! ={\text id}$. Therefore $f_*$ and
$g_*$ are both surjective and the exactness of (4.3) is equivalent to
the fact that

$${\ker} (g_* :A_p (Y') \rightarrow A_p (Y)) \ \overset {j_*}\to \to \
{\ker} (f_* :A_p (X') \rightarrow A_p (X))$$

\noindent is an isomorphism, hence to the fact that

$${\text {coker}} (A_p (Y) \overset {g^!}\to \to \ A_p (Y')) \
\overset {j_*} \to \to\ {\text {coker}} (A_p (X) \ \overset {f^*}\to
\to \ A_p (X'))$$

\noindent is an isomorphism. But this follows from the usual exact
sequence computing $A_p (X')$, \cite{F}, Proposition 6.7 c).

\medskip

To prove that (4.4) is exact we define once more $g^! :A^p (Y)
\rightarrow A^p (Y')$ by the formula $g^! (x)=g^* (x) \ c_{d-1} (F)$,
and we deduce from the equalities $f_* f^* ={\text id}$ ( \cite{F},
Proposition 17.5 a)) and $g_* \ g^! ={\text id}$ ( \cite{F} Example
3.3.3) and from the exact sequence

$$0 \rightarrow A^{p-d} (Y) \rightarrow A^{p-1} (Y') \oplus A^p (X)
\rightarrow A^p (X') \rightarrow 0$$

\noindent ( \cite{F} Example 17.5.1 c)) that we just need to show that
$j^* j_*$ induces an isomorphism

$${A^{p-1} (Y') \over g^! \ A^{p-d} (Y)} \rightarrow {A^p (Y') \over
g^* \ A^p (Y)} \ . \leqno (4.7)$$

\noindent Now $j^* j_*$ is the product by $\xi =c_1 ({\Cal O} (1))$ (
\cite{F} Proposition 17.4.1) and any element $x\in A^p (Y')$ can be
written uniquely as

$$x=\sum_{i=0}^{d-1} g^* (y_i ) \ \xi^i ,$$ \cite{F} Example 17.5.1
b), hence, also uniquely, as

$$x=\sum_{i=0}^{d-2} g^* (z_i ) \ \xi^i +g! (z),$$ since $c_{d-1}
(F)=\xi^{d-1} + \sum_{d-2}^{i=0} \ g^* (\alpha_i ) \
\xi^i$.

\noindent This implies that (4.7) is an isomorphism as in the proof of
i).

\bigskip

\subheading{ 4.5} To prove Theorem 5 iv) we let $U=X'-Y' =X-Y$ and we
consider the obvious morphism of exact sequences from $$ 0\to R^0 CH^p
(U) \to R^0 CH^p (X)\to R^0 CH^p (Y) \to R^1 CH^p (U) \to R^1 CH^p (X)
$$ to $$ 0\to R^0 CH^p (U) \to R^0 CH^p (X')\to R^0 CH^p (Y') \to R^1
CH^p (U)\to R^1 CH^p (X') $$ \noindent (see 3.1.1). Since $X$ is
smooth and projective we have $R^1 CH^p (X)=0$ (loc. cit.), therefore,
by diagram chase one concludes that (4.4) is exact. The same argument
using $L_i CH_p$ instead of $R^i CH^p$ proves that (4.3) is exact.
This ends the proof of Theorem 5.

\bigskip

\subheading{ 4.6} By applying Theorem 5, we can analyse the effect of
a blow up on the complex (3.1) considered in section 3.1.5.

Indeed, let $S$ be a regular noetherian scheme, $X$ a regular scheme
of finite type over $S$, $D\subset X$ a relative Cartier divisor with
normal crossings, and $Y_1 ,\ldots ,Y_n$ the components of the
associated reduced scheme $D^{\text{red}}$. As in 2.4 and 3.1.4 we may
define a cochain complex of abelian groups

$$C^.(X):A^p (X/S) \ \overset {\partial^*}\to \to \ A^p (Y^{(1)} /S) \
\overset {\partial^*}\to \to
 \ A^p (Y^{(2)}/S) \rightarrow \ldots$$

\noindent Now let $W\subset D^{{\text {red}}}$ be a regular
irreducible closed subscheme meeting the components $Y_i$ normally, in
the sense of \cite{H}, and let $f:X' \rightarrow X$ be the blow up of
$X$ along $W$.  The components of $f^* (D)^{{\text {red}}}$ are the
proper transforms $Y'_i$ of $Y_i$, $1\leq i\leq n$, together with the
exceptional divisor $Y'_0 =f^{-1} (Y)$. We define from these a complex
$C. (X')$ by the same method as $C. (X)$ (see 3.1.4).

\bigskip

\proclaim{ Proposition 6} The complexes $C^.(X)$ and $C^.(X')$ are
 quasi-isomorphic.  \endproclaim \bigskip

\demo{ Sketch of proof} We first define a new complex
 $B^.$ as follows.  If $I\subset \{ 0,\ldots ,n\}$ we let

$$\Omega_I =Y_I \quad \hbox{when} \quad 0\notin I$$

\noindent and

$$\Omega_I =W \cap Y_J =W_J \quad \hbox{if} \quad I=\{ 0\} \cup J.$$
Define

$$\Omega^{(r)} = \coprod_{{\text{card}} I=r} \Omega_I$$

\noindent and $B^r =A^p (\Omega^{(r)} /S)$. The differential
$\partial^*$ on $B^.$ are defined as for $C^.(X)$ from the inclusions
$\Omega_I \subset \Omega_J$ if $J\subset I$.

\medskip

There is a natural projection $B^. \rightarrow C^.(X)$.  Its kernel
$K^.$ is acyclic. Indeed we have $K^0 =0$ and $$ K^r =
\bigoplus_{{\text{card}} J=r-1}
 A^p (W_J/S)$$

\noindent when $r>0$. Let $P \subset \{ 1,\ldots ,n\}$ be the set of
indices $i$ such that $Y_i$ does not contain $W$. For any $I\subset P$
we get a subcomplex $K^._I$ of $K^.$ by considering only those $A^p
(W_J/S)$'s such that $J\cap P =I$. Since $W_J =W_I$ for all such $J$,
this complex $K^._I$ is acyclic. Using that $K^r = \ \oplus_{I\subset
P} K^r_I$ for all $r\geq 1$, we conclude that $K^.$ itself is acyclic
by applying the following lemma, the proof of which is left to the
reader: \bigskip

\proclaim{ Lemma 4} Let $A$ be a finite ordered set,
 let $C^.$ and $C^._{\alpha}, \alpha\in A$, be bounded cochain
complexes of abelian groups, and let $$f_r:\oplus_{\alpha \in
A}C^{r}_{\alpha}\to C^r$$ be group isomorphisms such that, for all
$x\in C^{r}_{\alpha}$, $(f_rd-df_r)(x)$ lies in $\oplus_{\alpha <
\beta}C^{r+1}_{\beta}$.  Assume that, for all $\alpha \in A$, the
complex $C^._{\alpha}$ is acyclic. Then $C^. $ is acyclic.
\endproclaim \bigskip

 On the other hand one checks that for any non empty subset $I\subset
\{ 1,\ldots ,n\}$ the intersection $Y'_I = \ \cap_{i\in I} \ Y'_i$ is
the blow up of $Y_I$ along $W_I =W\cap Y_I$, and that $Y'_0 \cap Y'_I$
is the exceptional divisor of that blow up (this uses the fact that
$W$ meets $D^{{\text {red}}}$ normally). In particular $f:X'
\rightarrow X$ maps $Y'_I$ onto $\Omega_I$ for any $I\subset \{
0,\ldots ,n\}$ and induces a map of complexes

$$f^* :B^. \rightarrow C^. (X').$$

\noindent One can prove that this is a quasi-isomorphism by
considering its cone $C(f^*)^.$ and the subcomplexes

$$0 \rightarrow A^p (Y_I /S) \rightarrow A^p (W_I /S) \oplus A^p (Y'_I
/S) \rightarrow A^p (Y_0 \cap Y'_I /S) \rightarrow 0$$

\noindent of $C(f^*)^.$, for all subsets
 $I\subset \{ 1,\ldots ,n\}$. We know from Theorem 5 iii) that these
complexes are acyclic, and, by applying Lemma 4 again, we conclude
that $f^*$ is a quasi-isomorphism.  \qed\enddemo

\bigskip

\subheading{ 4.7} One can use Proposition 6 to get another proof of
the results in 3.1.4. We know from \cite{H} that when $S= {\text
{Spec}} (k)$ and ${\text{ char}}(k)=0$, given two smooth
compactifications $X$ and $X'$ of a quasi-projective variety $U$ over
$k$, such that $X-U$ and $X'-U$ are divisors with normal crossing,
there exists a third one $X''$ and maps $X'' \rightarrow X$, $X''
\rightarrow X'$ which are the identity on $U$, one of them being the
composite of blow ups of the kind considered in Proposition 6. This
can be used to show that, up to quasi-isomorphism, $C. (X)$ depends
only on $U$.

\medskip To come back to motives, Theorem 5 implies that, given a
closed immersion $i : Y \rightarrow X$ of smooth proper varieties over
a field $k$, if $X'$ is the blow up of $X$ along $Y$ and $Y'$ its
exceptional divisor as in 4.1, the sequence of motives $$ 0
\rightarrow M(X) \overset {(i^* ,-f^*)}\to \rightarrow M(Y) \oplus
M(X') \overset {g^* + j^*}\to \rightarrow M(Y') \rightarrow 0 \leqno
(4.8) $$ is contractible. This follows from Theorem 5 iii) by Manin's
identity principle, as stated in \cite{Sc} 2.3 ii).

\medskip

This fact is the starting point of the alternative construction of
$W(X)$ due to Guillen and Navarro
\cite{G-N}, (5.4), my means of cubic hyperresolutions
\cite{G-N-P-P}.

\bigskip

\heading 5. $K$-Theory \endheading

\bigskip \subheading{5.1 Preliminaries}

\medskip

\subheading{5.1.1} We start by discussing the construction and some of
the properties of $K$-theory.  Recall that given an exact category
$\bold E$ in the sense of \cite{Q}, we can associate to it a category
$Q\bold E$, such that the classifying space $BQ\bold E$ .  Here by
``space'' we mean fibrant simplicial set, and for a category $\bold
C$, $B\bold C$ denotes the result of applying Kan's ${\text
{Ex}}^\infty$ functor to the nerve $N.\bold C$ of $\bold C$.  We shall
assume that the zero objects of all exact categories that we deal with
are unique. (If necessary, given an arbitrary exact category $\bold
E$, we may form an equivalent exact category $\tilde{\bold E}$ with a
unique zero object by replacing the subcategory of zero objects of
$\bold E$ by a single zero object.)  The classifying space $BQ\bold E$
is therefore canonically pointed. We furthermore assume that this
space is the zeroth space of a spectrum $\bold K(\bold E)$, in the
sense of Appendix A.1, with the following properties:

\item{$\bullet$} ${\bold E}\mapsto {\bold K(\bold E)}$
 is strictly (i.e. not just up to homotopy) functorial;
\item{$\bullet$} Given exact categories $\bold E$, $\bold F$, and
$\bold G$, and a biexact functor $$\mu : {\bold E}\times {\bold F}\to
{\bold G}$$ there is a canonical, functorial, pairing: $${\bold
K}(\mu) : {\bold K(\bold E)}
\wedge {\bold K(\bold F)}\to {\bold K(\bold G)}\quad .$$

\noindent
For example we can use the multiple $Q$-construction; see \cite{G1}.
In particular, as in \cite{T1} Appendix A, we shall assume that these
spectra are fibrant, and that they are cofibrant as prespectra (i.e.
the map of the suspension of the $i$-th space to the $(i+1)$-st space
is injective).  The K-theory groups are then defined as the homotopy
groups of these spectra: $$K_m({\bold E})=\pi_{m+1}({\bold K(\bold
E)}),$$ and are functorial with respect to exact functors between
exact categories.  They are also compatible with products in the sense
that, given a biexact functor $$\mu : {\bold E}\times {\bold F}\to
{\bold G}$$ as above, we get a functorial pairing of graded groups
$$K_*(\mu):K_m({\bold E})\otimes K_n({\bold F})\to K_{m+n}({\bold
G})\quad .$$

\bigskip \subheading{5.1.2} Given a (noetherian) scheme $X$, we can
consider two exact categories: the category ${\bold M}(X)$ of coherent
sheaves of ${\Cal O}_X$-modules and the sub-category ${\bold P}(X)$ of
${\bold M}(X)$ consisting of locally free modules. We then obtain the
groups: $$K_m(X)=K_m({\bold P}(X))$$ and $$K'_m(X)=K_m({\bold
M}(X))\quad .$$

 These groups are functorial with respect to pull-back and proper
push-forward respectively.  However the $K$-theory spectra themselves
are not strictly functorial, but rather are functorial only up to
homotopy, and similarily the projection formula is also only true up
to homotopy; this is because the underlying functors from the category
of varieties to the category of categories are ``lax'' rather than
``strict'' . To remedy to this we must rigidify the underlying
category valued functors.  The standard constructions for rigidifying
lax functors are due to Street, \cite{S}. We now give a description of
Street's (second) construction in the case of locally free sheaves.
We then describe, following Thomason \cite{T-T}, an intrinsically
rigid construction of $K'$-theory; we then modify this construction in
order to make the projection formula true exactly rather than up to a
natural isomorphism.

 Given a scheme we let ${\bold P}^{\text{Big}}(X)$ be the category of
locally free sheaves in the big Zariski site over $X.$.  An object in
${\bold P}^{\text{Big}}(X)$ consists of a locally free sheaf ${\Cal
F}_f$ on $Y$ for each map of schemes $f:Y\to X$, and of an isomorphism
$g^*({\Cal F}_{f_2})\to {\Cal F}_{f_1}$ for each morphism
$g:(f_1:Y_1\to X)\to (f_2:Y_2\to X)$ in the category of schemes over
$X$ (i.e. a map $g:Y_1\to Y_2$ such that $f_2\circ g= f_1$), with the
obvious compatibility with respect to composition.  We omit the proof
of the following proposition, since it is straightforward.

\proclaim{Lemma 5} The forgetful functor from ${\bold
P}^{\text{Big}}(X)$ to the category of locally free sheaves on (the
small Zariski site of) $X$ is an equivalence of categories.
\endproclaim

 If $f:X\to Y$ is a map of schemes, then we have a restriction functor
from $f^*:{\bold P}^{\text{Big}}(Y) \to {\bold P}^{\text{Big}}(X)$
which takes a family $(g:Z\to Y)\to{\Cal F}_g$ to the family $(h:Z\to
X)\to {\Cal F}_{h\circ f}$.  Clearly if $g:V\to X$ is another map,
then $(fg)^*=g^*f^*$ is an equality of functors, not just a natural
equivalence. Under the previous equivalence of categories this functor
is compatible with the usual pull-back map on vector bundles.
Therefore if we define ${\bold K}(X)$ to be the spectrum ${\bold
K}({\bold P}^{\text{Big}}(X))$ we obtain a (strictly) contravariant
functor from schemes to spectra such that $$K_m(X)= \pi_{m+1} {\bold
K}({\bold P}^{\text{Big}}(X))$$ for all $m\geq 0$.

\bigskip \subheading{5.1.3} In order to make $K'$
 covariant functorial, we use the
construction of Thomason
\cite{T-T}. Given a scheme $X$ let
 ${\bold C}^b(X)$ be the category of complexes of
flasque quasi-coherent sheaves of ${\Cal O}_X$-modules, having
cohomology that is coherent and bounded.  Taking the category $w$ of
quasi-isomorphisms of sheaves to be the weak equivalences, and the
standard notion of exact sequence, the pair $({\bold C}^b(X),w)$ is a
category with cofibrations and weak equivalences.  The following is
due to Thomason, {op. cit}:

\proclaim{Proposition 7}

\item{i)} $K_*({\bold C}^b(X),w)\simeq K'(X)$
\item{ii)} If $f:X\to Y$ is a proper morphism
of schemes, ${\Cal A}\mapsto f_*{\Cal A}$
is an exact functor, preserving weak equivalences, $f_*: ({\bold
C}^b(X),w) \to ({\bold C}^b(Y),w)$.  Furthermore, if $g: Y\to Z$, then
we have an identity (not just a natural equivalence) of functors
$g_*f_*=(gf)_*$

\endproclaim

 To make the projection formula itself, and not just its constituent
functors, an identity, we must still rigidify further.  We use a
construction similar, but not identical, to Street's first
construction in \cite{S}.  Let $\tilde{\bold C}(X)$ denote the
category with objects pairs $(f:Y\to X,{\Cal A}^.\in {\bold
C}^b(Y))=(f,{\Cal A}^.)$ and morphisms $$\theta\in {\text
{Hom}}_{\tilde{\bold C}(X)}((f,{\Cal A}^.),(g,{\Cal B}^.))$$ given
simply by maps $\theta':f_*{\Cal A^.}\to g_*{\Cal B^.}$of complexes of
sheaves on $X$.  We say that $\theta$ is a weak equivalence if
${\theta}'$ is a quasi-isomorphism.  Similarily we say that $\theta$
is a cofibration if ${\theta}'$ is a monomorphism equal to the kernel
of a map of complexes in ${\bold C}^b(X)$.

\proclaim{Lemma 6} The obvious inclusion functor
$j:{\bold C}^b(X)\to \tilde{\bold C}(X)$ is an equivalence of
categories with cofibrations and weak equivalences.  In particular, it
induces a homotopy equivalence of $K$-theory spectra.

\endproclaim

\demo{Proof}
Clearly the inclusion functor preserves cofibrations and weak
equivalences.  Now consider the functor $p:\tilde{\bold C}^b(X)\to
{\bold C}(X)$ given by $$p:(f,{\Cal A}^.)\to f_*{\Cal A}^.$$ on
objects, and by the identity on Hom-sets.  The composition $p\circ j$
is clearly the identity. In the other direction, the composition
$j\circ p$ is isomorphic to the identity functor on $\tilde{\bold
C}^b(X)$, mapping $(1_X,f_*{\Cal A^.})$ to $(f, \Cal A^.)$.

 \qed \enddemo

 Now, given a map of schemes $g:X\to Y$, we have an exact functor
$g_*:\tilde{\bold C}^b(X)\to \tilde{\bold C}^b(Y)$, given on objects
by: $$g_*:(f,{\Cal A^.})\mapsto (g\circ f, {\Cal A^.})$$ and on
morphisms by the natural action of $g_*$ on the underlying Hom-sets of
complexes on $X$.  One can check that, via the equivalence of
categories of Lemma 6, this functor is compatible with the usual
direct image.

 The usual cap-product: $K_m(X)\otimes K'_n(X)\to K'_{m+n}(X)$ is
usually viewed as being induced by the bi-exact functor $$
\otimes_{{\Cal O}_X}:{\bold P}(X)\times {\bold M}(X)\to {\bold M}(X).
$$ Given a proper morphism $f:X\to Y$, we have the projection formula,
for $\alpha\in K_*(Y)$, and $\beta\in K'_*(X)$: $$ \alpha\cap
f_*(\beta)=f_*(f^*(\alpha))\cap \beta) \quad .  $$ This formula is
usually derived from the isomorphism of functors: $$ {\Cal
F}\otimes_{{\Cal O}_Y}f_*{\Cal M}\simeq f_*({\Cal F}\otimes_{{\Cal
O}_X}{\Cal M})
\quad .
$$ Now we can describe the cap-product as follows: there is a biexact
functor $$\mu:{\bold P}^{\text{Big}}(X)\times {\bold C}^b(X) \to
{\bold C}^b(X)$$ $$\mu:((g:Z\to X)\to{\Cal F}_g)\times (f:Y\to X,{\Cal
A}^.)
\mapsto
(f:Y\to X,{\Cal F}_f\otimes {\Cal A}^.)  .$$ On morphisms this acts as
follows. A morphism of locally free sheaves in the big Zariski site
$(f\mapsto \phi_f) : (f\mapsto{\Cal F}_f)\to(f\mapsto {\Cal G}_f)$
induces the morphism $$(\phi_f\otimes 1):(f:Y\to X,{\Cal F}_f\otimes
{\Cal A}^.) \to (f:Y\to X,{\Cal G}_f\otimes {\Cal A}^.)  $$ in ${\bold
C}^b(X)$.  A morphism in the category $ {\bold C}^b(X)$
$$\theta:(f,{\Cal A}^.)\to (g,{\Cal B}^.),$$ i.e. a morphism of
complexes of sheaves on $X$ $$\theta':f_*{\Cal A}^.\to g_*{\Cal B}^.$$
gives a map $$(1,\theta'):((h\mapsto{\Cal F}_h),(f,{\Cal A}^.))\to
((h\mapsto{\Cal F}_h),(g,{\Cal B}^.))$$ in the product category
${\bold P}^{\text{Big}}(X)\times {\bold C}^b(X)$ and its image by the
functor $\mu$ is the unique map: $$f_*({\Cal F}_f\otimes {\Cal
A}^.)\to g_*({\Cal F}_g\otimes {\Cal B}^.)$$ which makes the diagram
$$
\CD
f_*({\Cal F}_f\otimes {\Cal A}^.)@>>> g_*({\Cal F}_g\otimes {\Cal
B}^.)\\ @VVV @VVV\\ {\Cal F}_{id_X}\otimes f_*{\Cal A}^.@>>\theta >
{\Cal F}_{id_X}\otimes g_*{\Cal B}^.
\endCD
$$ commutative, where the vertical maps are the isomorphisms induced
by the projection formula at the level of modules.

 Given a map of varieties $h:X\to Y$, an object $A=(g\mapsto {\Cal
F_g}\in {\bold P}^{\text Big}(Y))$ and an object $B=(f:Z\to Y,{\Cal
A}^.)$ in ${\bold C}^b(Y)$ we find that : $$
h_*(\mu(h^*A,B))=h_*(f:Z\to Y,{\Cal F}_{hf}\otimes {\Cal A}^.)  =
(h\cdot f :Z\to Y,{\Cal F}_{hf}\otimes {\Cal A}^.)=\mu (A,h_*B) $$ One
may also check compatibility for morphisms.

To summarize, we have functors $K$ and $K'$ from the category of
projective varieties to ${\bold {Spectra}}$, the first contravariant,
the second covariant, together with products which satisfy the
projection formula exactly rather than up to homotopy.  For example,
if $X.:\Delta^{op}\to {\bold {Pr}}$ is a simplicial projective
variety, we get from $X.$ by the construction above a simplicial
spectrum ${\bold K'}:\Delta^{op}\to {\bold {Spectra}}$ and we can
define the $K'$-prespectrum of $X.$ as being the corresponding
homotopy colimit: $${\bold K'}(X.) = {\text {hocolim}}_
{\Delta^{op}}{\bold K'}(X_p).$$ Its homotopy groups $K'_m(X.) =
\pi_{m+1}{\bold K'}(X.)$ are the abutment of a first-quadrant
convergent spectral sequence $$E_{pq}^2 = H_p \left(*\mapsto K'_q (X_*
)\right)
\Rightarrow K'_{p+q}(X.)$$ (\cite{B-K} XII 5.7 and \cite{T1}
Proposition 5.7).

 When dealing with the $K$-theory of simplicial schemes one gets a
similar definition by replacing homotopy colimits with homotopy
limits.  However the associated spectral sequences need not be
convergent in general (see \cite{B-K} XII 7 or \cite{T1} 5.44).  The
consideration of $K_0$-motives in the next paragraph will help us to
solve this difficulty in Section 5.3.

\bigskip \subheading{5.1.4
Remark} When $X$ is smooth $K_m(X)=K'_m(X)$ for all $m\geq 0$
\cite{Q}. In \cite{G1} Lemma 4.5 it is asserted that in general any
element in $K_m(X)$ is the inverse image of an element in $K_m(M)$,
where $M$ is a smooth variety.  This is used in \cite{So2} 6.2 to
define operations on the $K$-theory of singular varieties.  However
the proof of \cite{G1} Lemma 4.5 is incorrect since the compatibilty
statements (c) in loc.cit., p.247, are not enough to describe an
arbitrary diagram in $Q{\bold P}(X)$.

\bigskip \subheading{5.2 $K_0$-motives} \medskip
\subheading{5.2.1}
 We shall now give analogs of the results in Section 1 for $K$-theory
instead of Chow groups. From now on ${\bold V}$ will denote the
category of smooth projective varieties over a fixed field $k$.

 First remark that one gets a theory of motives by replacing Chow
groups by $K_0$ in the definition of correspondences (see also
\cite{M1}). Namely, let
 $\bold {KC}$ be the category with the same objects as ${\bold V}$,
and with morphisms $$ {\text {Hom}}_{{\bold {KC}}} (X,Y) = K_0 (X
\times Y).  $$ The composition law is defined as in ${\text
{Hom}}_{{\bold C}} (X,Y)$ (see 1.2) and there is a covariant functor
${\bold V} \rightarrow \bold {KC}$ mapping a morphism $f: X
\rightarrow Y$ to the class $[ {\Cal O}_{\Gamma_f}] \in K_0 (X \times
Y)$ of the structure sheaf of the graph of $f$. Let ${\bold {KM}}$ be
the associated category of motives, defined as in 1.3; we call them
{\it $K_0$-motives}.

\medskip

 Notice that for all $m\geq 0$, the functor $K_m$ from ${\bold V}$ to
abelian groups can be factored through ${\bold {KC}}$, and hence
$\bold {KM}$, both as a covariant and a contravariant functor.

\medskip

 We obtain results similar to those in Section 1 by replacing the
Gersten complexes $R_{q,*}$, $q\geq 0$, with the $K'$-theory spectrum.
However, since the functor $X \mapsto {\bold K'}(X)$ does not factor
through ${\Bbb Z} {\bold V}$ we need to modify some of the arguments.
Instead of complexes in ${\Bbb Z} {\bold V}$ we must work with
simplicial schemes. The following will play the role of Theorem 1:

\bigskip

\proclaim{Theorem 6} Let $$ \CD X.@>>> Y.\\
    @VVV @VVV \\ S. @>>> T.  \endCD $$ be a commutative square of maps
between simplicial objects in ${\bold V}$. Suppose that, for all
varieties $V$ in ${\bold V}$, the associated square of spectra $$ \CD
{\bold K}'(V\times X.)@>>> {\bold K}'(V\times Y.)\\ @VVV @VVV \\
{\bold K}'(V\times S.) @>>> {\bold K}'(V\times T.)
\endCD \leqno (5.1) $$ is homotopy cartesian. Then the associated
square of complexes of $K_0$-motives $$ \CD KM(S.) @<<< KM(T.)\\ @VVV
@VVV \\ KM(X.)@<<< KM(Y.)\endCD $$ is homotopy cartesian (i.e. the
associated total complex is contractible).

\endproclaim

\bigskip

\demo{ Proof } The square $(5.1)$ is homotopy cartesian if and only if
the homotopy colimit of the diagram of simplicial spectra $$ \CD
{\bold K}'(i\mapsto V\times X_i)@>>> {\bold K}'(i\mapsto V\times
Y_i)@>>> *\\ @VVV @VVV @.\\ {\bold K}'(i\mapsto V\times S_i)@>>>
{\bold K}'(i\mapsto V\times T_i)@.\\ @VVV @.  @.\\ *@.@.\\ \endCD $$
is contractible. By \cite{B-K} XII 3.3., this iterated homotopy
colimit is isomorphic to the homotopy colimit of the associated $I
\times \Delta^{\text {op}}$-diagram of spectra, where $I$ is the
finite category $$ \CD a@>>> b@>>> *\\ @VVV @VVV @.\\ c@>>> d\\ @VVV
@. \\ *@.@.\ .\\ \endCD $$ If $A.$ denotes this diagram, there is an
associated spectral sequence $$ E_{pq}^2 = H_p \left( I
\times \Delta^{\text {op}} , \pi_q ( A.)  \right) \Rightarrow
\pi_{p+q} \ {\text {hocolim}}_{I \times
\Delta^{\text {op}}} (A.)  $$ ( \cite{B-K} XII 5.7 and \cite{T1}
Proposition 5.17). Given any functor $\Phi$ from $I \times
\Delta^{\text {op}}$ to abelian groups which vanishes on $* \times
\Delta^{\text {op}}$, the homology groups
 $H_* (I \times \Delta^{\text {op}}, \Phi)$ are those of the homotopy
push-out of the diagram of chain complexes $$ \CD \Phi(a,.)@>>>
\Phi(b,.)@>>> 0\\ @VVV @VVV @.\\ \Phi(c,.)@>>> \Phi(d,.)@.\\ @VVV @.
@.\\ 0@.@.\ ,\\
\endCD $$
 i.e. the total complex of the square $$ \CD \Phi(a,.)@>>> \Phi(b,.)\\
@VVV @VVV \\ \Phi(c,.)@>>> \Phi(d,.)\ .\\ \endCD $$ Thus the groups
$E_{pq}^2$, $p\geq 0$, are the homology groups of the complex $K_q (V
\times C. )$, where $C.$ is the total complex in ${\Bbb Z} {\bold V}$
associated to the commutative square of simplicial varieties $$ \CD
X.@>>>Y.\\ @VVV @VVV \\ S. @>>> T.\ .
\endCD
$$ It follows that, for all $V$, we have a spectral sequence $$
E_{pq}^2 = H_p \left( K_q (V \times C_* )\right) \Rightarrow 0 .  $$
As in the proof of Theorem 1, we may now prove by induction on $n\geq
1$
that
$C.$ is contractible as a complex of $K_0$-motives in degrees less
than $n$, i.e. that there exist $K_0$-correspondences $h_i \in {\text
{Hom}}_{{\bold {KC}}} (C_i ,C_{i+1})$, such that $$ h_{i-1} \circ
\delta_i +\delta_{i+1} \circ h_i = 1_{C_i} \ , \ 0 \leq i \leq n-1 .
$$ Since $K_q$ factors via ${\bold {KC}}$ we get that $E_{pq}^2 = 0$
for all $q\geq 0$ and $p\leq n-1$. The end of the proof is then
parallel to that of Theorem 1.
\qed \enddemo \bigskip

\subheading{5.2.2}
 Using Theorem 6 and the descent Theorem 4.1 in \cite{G2} for
$K'$-theory, we can associate to any variety $X$ over a field $k$ of
characteristic zero a cochain complex of $K_0$-motives $KW(X)$ in the
homotopy category ${\text {Hot}} ({\bold {KM}})$, which is
well-defined up to canonical isomorphism and enjoys the same
properties as $W(X)$ in Theorem 2 above. If $\overline X$ is a
compactification of $X$ and $\widetilde j : \widetilde{Y}. \rightarrow
\widetilde{X}.$ a non-singular hyperenvelope of the inclusion $j :
\overline X - X
\rightarrow \overline X$, $KW(X)$ is represented by the complex $$ C
\left( KM \left( \widetilde{X}. \right) \overset {\widetilde{j}^*}\to
\rightarrow KM
\left( \widetilde{Y}. \right) \right)^. [-1] \quad
\hbox{in} \quad {\text {Hot}} ({\bold {KM}}).  $$

 The proof of the properties of $KW(X)$ is the same as in Theorem 2.
For instance let $\widetilde j \rightarrow \widetilde{j}'$ be a map of
non-singular hyperenvelopes of $j$. To check that these define the
same complex $KW(X)$ up to homotopy equivalence, notice that for any
$V$ in ${\bold V}$ the associated square of $K'$-theory spectra
${\bold {K'}} \left( 1_V \times \widetilde j \right) \rightarrow
{\bold {K'}} \left( 1_V \times j
\right)$, ${\bold {K'}} \left( 1_V
\times \widetilde{j}' \right) \rightarrow {\bold {K'}} \left( 1_V
\times j \right)$ and hence ${\bold {K'}} \left( 1_V
\times \widetilde{j}' \right)
 \rightarrow {\bold {K'}} \left( 1_V \times \widetilde j \right)$ are
homotopy cartesian. Therefore, by Theorem 6, the associated square of
complexes of $K_0$-motives $KM\left( \widetilde j \right) \rightarrow
KM \left( \widetilde{j}' \right)$ is homotopy cartesian, as desired.

 Similarly, the same proof as in 2.4 tells us that $KW(X)$ is
represented by a bounded complex in ${\bold {KM}}$, of length at most
$\dim (X)+1$.

\bigskip

\subheading{5.3 $K$-theory with compact support}

\medskip \subheading{5.3.1}
 We are now able to define the {\it $K$-theory with compact support}
of any variety over a field $k$ of characteristic zero.

\medskip

 First consider the case of a complete variety $X$ over $k$. Let
$\widetilde{X}. \rightarrow X$ be a non-singular hyperenvelope of $X$.
By applying the contravariant $K$-theory functor ${\bold K} : {\bold
V} \rightarrow {\bold {Spectra}}$ we get a cosimplicial spectrum
Following \cite{T1}, section 5.6, we can form the spectrum ${\bold
K}(X.):={\text {holim}}_n \ K(\tilde{X}_n)$.  By op. cit.  Proposition
5.13, the spectral sequences of Bousfield and Kan in the unstable case
(\cite{B-K} IX 5 and XII 7.1) give rise to an unfringed spectral
sequence abutting to $\pi_*({\bold K}(X.))$

\medskip

 We claim that this spectral sequence is strongly convergent. Indeed
its $E_2$-term is $$ E_2^{pq} = H^p \left(*\mapsto \pi_{-q} {\bold K}
\left(\widetilde{X}_* \right) \right) = H^p \left(*\mapsto K_{-q}
\left(\widetilde{X}_* \right)
\right) ,
$$ and we know from the last section that the complex of $K_0$-motives
$KM \left( \widetilde{X}. \right)$ is homotopy equivalent to a bounded
complex of length at most $\dim (X) +1$. Therefore $E_2^{pq} = 0$ when
$p\geq\dim (X)+1$, and, for all $r\geq 2$, $E_r^{pq} = 0$ unless $q+i
\leq 0$, $0\leq p
\leq \dim
 (X)$, and $p+q \leq 0$. The spectral sequence therefore converges
strongly to $\pi_{-p-q} \ {\text {holim}}_n \ K \left(
\widetilde{X}_n \right)_i$ (see \cite{B-K} loc. cit. or \cite{T1}
Lemma 5.48 i)).

 In general, we get a functor from simplicial schemes to the category
of spectra $$X.\mapsto {\bold K}(X.):={\text {holim}}_n \ {\bold
K}(X_n)\quad .$$ We then define $K_*(X.):=\pi_*({\bold K}(X.))$.  Note
that for $\tilde{X}.$ a hyperenvelope of a variety $X$, the spectral
sequence above will be concentrated in the strip $0 \leq p \leq \dim
(X)$ and that therefore $K_m \left( \widetilde{X}.
\right)$ will in general be non-zero for negative values of $m$,
$m\geq -\dim (X)$.

\medskip

 Now let $X$ be an arbitrary variety over $k$, $X \subset \overline X$
a compactification, $j:Y=\overline X -X \hookrightarrow \overline X$
its complement, and $\widetilde j : \widetilde{Y}. \rightarrow
\widetilde{X}.$ a non-singular hyperenvelope of $j$. Consider the
homotopy fiber ${\bold K} \left( \widetilde j \right)$ of the map of
spectra $\widetilde{j}^* : {\bold K} \left( \widetilde{X}. \right)
\rightarrow {\bold K} \left( \widetilde{Y}. \right)$. We view it as
the homotopy limit of the diagram $$ \CD *@.\\ @VVV @.\\ {\bold
K}(\widetilde{Y}.) @<<< {\bold K}(\widetilde{X}.)\ , \endCD $$ which
is isomorphic by
\cite{B-K} XI 4.3.  to the homotopy limit of the
 corresponding $I \times \Delta$ diagram of spectra, where $I$ is the
small category $$ \CD *@.\\ @VVV @.\\ a @<<<b\ .  \endCD $$ Therefore
we get a convergent spectral sequence $$ H^p \left(*\mapsto C
\left( K_{-q} \left( \widetilde{X}.\right) \rightarrow K_{-q} \left(
\widetilde{Y}. \right) \right)^{*-1} \right) \Rightarrow K_{-p-q}
\left( \widetilde j \right) , \leqno (5.2) $$ where $K_m \left(
\widetilde j \right) = \pi_m \ {\bold K} \left( \widetilde j \right)$.

\medskip

 The spectrum ${\bold K} \left( \widetilde j \right)$ is our
definition of the $K$-theory with compact support of $X$. To see that
is independent of choices up to canonical homotopy equivalence,
consider once more two compactifications $\overline{X}_1$ and
$\overline{X}_2$ of $X$, $j_1 :
\overline{X}_1 -X \rightarrow \overline{X}_1$ and $j_2
 : \overline{X}_2 -X \rightarrow \overline{X}_2$ their complements,
$\pi : \overline{X}_1
\rightarrow \overline{X}_2$ a morphism which is the
 identity on $X$, and $\widetilde{\pi} : \widetilde{j}_1 \rightarrow
\widetilde{j}_2$ a map of non-singular hyperenvelopes of $j_1$ and
$j_2$, compatible with the morphism $\pi : j_1 \rightarrow j_2$ in
${\text {Ar}} ({\bold P})$ in the obvious way (compare 2.2 and 2.3
above). There is then a map of spectral sequences from $$ H^p
\left(*\mapsto C
\left( K_{-q} \left( \widetilde{X}_{2}. \right) \rightarrow K_{-q}
\left( \widetilde{Y}_{2}. \right) \right)^{*-1} \right) \Rightarrow
K_{-p-q} \left( \widetilde{j}_2 \right) $$ to $$ H^p \left( *\mapsto C
\left( K_{-q} \left( \widetilde{X}_{1}. \right) \rightarrow K_{-q}
\left(
\widetilde{Y}_{1}. \right) \right)^{*-1} \right) \Rightarrow K_{-p-q}
\left( \widetilde{j}_1 \right) $$ induced by $\widetilde{\pi}$. Since
the map of complexes of motives $$ C \left( KM \left(
\widetilde{X}_{2}. \right) \rightarrow KM \left( \widetilde{Y}_{2}.
\right)\right)^. \rightarrow C \left( KM \left( \widetilde{X}_{1}.
\right) \rightarrow KM
 \left( \widetilde{Y}_{1}. \right)\right)^.  $$ is a homotopy
equivalence, the map of $E_2$-terms in the above spectral sequences is
an isomorphism and therefore $\widetilde{\pi}^* : {\bold K} \left(
\widetilde{j}_2 \right) \rightarrow {\bold K}
\left( \widetilde{j}_1 \right)$ is a homotopy equivalence.

\medskip

 As in Section 2.3 we can then show that, given any proper map $f:X_1
\rightarrow X_2$ of varieties, together with compactifications $X_i
\hookrightarrow \overline{X}_i$, $j_i : \overline{X}_i -X_i
\rightarrow \widetilde{X}_i$, and non-singular hyperenvelopes $\pi_i :
\widetilde{j}_i \rightarrow j_i$, there is a canonical map in the
homotopy category of spectra $$ {\bold K} (f) : {\bold K} \left(
\widetilde{j}_2 \right) \rightarrow {\bold K} \left( \widetilde{j}_1
\right) , $$ which is a homotopy equivalence when $f$ is an
isomorphism. If $f :X_1 \rightarrow X_2$ and $g:X_2 \rightarrow X_3$
are two such maps, and if we choose compactifications and
hyperenvelopes $\widetilde{j}_i$, $i=1,2,3$, for all three varieties,
then $$ {\bold K} (gf) = {\bold K} (f) \ {\bold K} (g) : {\bold K}
\left( \widetilde{j}_3 \right) \rightarrow {\bold K} \left(
\widetilde{j}_1 \right) .  $$ It follows that, if we choose a
compactification and a non-singular hyperenvelope $\widetilde{j}_X$ of
its complement for every variety $X$, we obtain a contravariant
functor ${\bold K}^c$ from varieties to the homotopy category of
spectra (i.e. the stable homotopy category) by sending $X$ to ${\bold
K} \left( \widetilde{j}_X
\right)$ and any proper morphism $f$
 to ${\bold K} (f)$. Two families of choices give rise to canonically
isomorphic functors.

\medskip

 The properties of this functor ${\bold K}^c$ are summarized in the
following theorem : \bigskip

\proclaim{Theorem 7} To each variety $X$ over a field $k$ of
characteristic zero is associated a spectrum ${\bold K}^c (X)$, which
is well defined up to canonical homotopy equivalence and enjoys the
following properties: \smallskip

\item{i)} {\it If $X$ is complete and non-singular, ${\bold K}^c (X)$
is the usual Quillen $K$-theory of vector bundles.}

\smallskip
\item{ii)} {\it Any proper map $f:X \rightarrow X'$ of varieties
induces a pull-back map of spectra $f^* : {\bold K}^c (X') \rightarrow
{\bold K}^c (X)$.  Given two composable proper maps $f$ and $g$, then
$(fg)^* = g^* f^*$.}
\smallskip
\item{iii)} {\it Any open immersion $i:U \rightarrow X$ induces a map
of spectra $i_* : {\bold K}^c (U) \rightarrow {\bold K}^c (X)$. Given
two composable open immersions $i$ and $j$, then $(ji)_* = j_* i_*$.}
\smallskip
\item{iv)} {\it If $i:U\rightarrow X$ is an open immersion with
complement $j:Y=X-U \rightarrow X$, there is a fibration sequence of
spectra} $$ {\bold K}^c (U) \overset {i_*}\to \rightarrow {\bold K}^c
(X) \overset {j^*}\to \rightarrow {\bold K}^c (Y).  $$ \smallskip
\item{v)} {\it If $K_m^c (X) = \pi_m \ {\bold K}^c (X)$ denote the
homotopy groups of ${\bold K} (X)$ then $K_m^c (X) = 0$ if $m<-\dim
(X)$. There is a strongly convergent weight spectral sequence $$
E_r^{pq} (X)
\Rightarrow K_{-p-q}^c (X) \ , \quad r\geq 2 , $$ which
 is equal to the spectral sequence (5.2) above for any choice of a
compactification of $X$ and of a non-singular hyperenvelope
$\widetilde j$. The associated filtration $F^pK_m^c (X)$, called the
{\rm {weight filtration}}, is increasing, finite, and independent of
choices.

\smallskip

\item{vi)} If ${\Bbb A}^1$ is the affine line on $k$, the inclusion of
$X$ as $X \times \{ 0 \}$ in $X \times {\Bbb A}^1$ induces
isomorphisms $K_m^c (X) \simeq K_m^c (X\times {\Bbb A}^1)$.
\endproclaim \bigskip

\subheading{5.3.2}
Finally, we shall describe pairings of spectra between $K$-theory with
compact support and $K'$-theory of varieties.

\medskip

 Recall from 5.1 that ${\bold K}'$ and ${\bold K}$ are covariant and
contravariant functors from the category of projective varieties to
the category of spectra and that the cap-product $$ {\bold K} (X)
\wedge {\bold K}' (X) \rightarrow {\bold K}' (X) $$ satisfies the
projection formula exactly (and not up to homotopy). Therefore, if
$X.$ is a simplicial object in ${\bold V}$, there is a pairing of
$\Delta^{op}$-diagrams $$ {\bold K} (X.) \wedge {\bold K}' (X.)
\rightarrow {\bold K}' (X.)  $$ in the sense of the Appendix A.2.1,
and hence, by Proposition 9, a pairing of prespectra $$ {\text
{holim}}_n
\ {\bold K} (X_n) \wedge {\text {hocolim}}_n \ {\bold K}' (X_n)
\rightarrow {\text {hocolim}}_n \ {\bold K}' (X_n) .  $$ For any map
$f:X. \rightarrow Y.$ of simplicial objects, this pairing satisfies
the projection formula exactly.

\medskip

 Now let $X$ be a variety over a field $k$ of characteristic zero,
$\overline X$ a compactification of $X$, $j:\overline X - X
\rightarrow \overline X$ the complementary inclusion, and $\widetilde
j : \widetilde{Y}.  \rightarrow \widetilde{X}.$ a non-singular
hyperenvelope of $j$. We obtain a commutative diagram of pairings of
spectra $$\CD {\bold K} \left( \widetilde{X}.  \right)\wedge {\bold
K}' \left(
\widetilde{Y}.  \right) @>{\widetilde{j}^*\wedge \text{Id}}>> {\bold
K} \left( \widetilde{Y}.  \right)\wedge {\bold K}' \left(
\widetilde{Y}.  \right) @>\cap>> {\bold K}' \left( \widetilde{Y}.
\right) \\
      @VV\text{Id}\wedge \widetilde{j}_* V @.  @VV\widetilde{j}_* V \\
{\bold K} \left( \widetilde{X}.  \right)\wedge {\bold K}' \left(
\widetilde{X}.  \right)
  @>>> @>>> {\bold K}' \left( \widetilde{X}.  \right) \ .  \endCD $$
Therefore, by the Appendix A.2.2, we have a pairing $$ {\text {Fiber}}
\left( \widetilde{j}^* \right) \wedge {\text {Cofiber}}
\left( \widetilde{j}_* \right) \rightarrow {\bold K}' \left(
\widetilde{X}.  \right) \overset {\sim} \to \rightarrow{\bold K}'
\left( \overline X \right) , $$ and hence a pairing of spectra $$ \mu
: {\bold K}^c (X) \wedge {\bold K}' (X) \rightarrow {\bold K}' \left(
\overline X \right) $$ for which one can check the following
projection formulae:

\bigskip

\proclaim{Proposition 8} Assume $f:X \rightarrow Y$ is a proper map of
varieties, and $\overline f : \overline X \rightarrow \overline Y$
extends $f$ to compactifications of $X$ and $Y$. Then the following
diagram is commutative $$\CD {\bold K}^c (Y)\wedge {\bold K}' (X) @>>>
{\bold K}^c (X)\wedge {\bold K}' (X) @>>> {\bold K}' \left( \overline
X \right) \\ @VVV @.  @VVV \\ {\bold K}^c (Y)\wedge {\bold K}' (Y)
@>>>\longrightarrow @>>> {\bold K}' \left( \overline Y \right)\ .  \\
\endCD $$

\bigskip
 If $i:U \rightarrow X$ is an open immersion and $\overline X$ a
compactification of $X$, the following diagram is commutative $$\CD
{\bold K}^c (U)\wedge {\bold K}' (X) @>>> {\bold K}^c (X)\wedge {\bold
K}' (X) \\ @VVV @VVV \\ {\bold K}^c (U)\wedge {\bold K}' (U) @>>>
{\bold K}' \left( \overline X \right)\ .  \\ \endCD $$ \endproclaim
Finally, let us remark that we can define the weight filtration
$F_qK'_n(X)$ as the filtration coming from the spectral sequence $$H_q
\left(*\mapsto C \left( K'_{r} \left( \widetilde{Y}. \right) \rightarrow
K'_{r} \left( \widetilde{X}. \right) \right)_* \right)
\Rightarrow K'_{q+r} (X)$$ as in (5.2) using the descent theorem for
$K'$-theory \cite{G2} Theorem 4.1.  This descending filtration is
finite and independent of choices.  We expect the pairing of spectra
considered above to induce a pairing on the weight filtrations:

$$F^pK^c_m(X)\otimes F_qK'_n(X)\to F_{q-p}K'_{m+n}(\overline X).$$

\bigskip \heading Appendix: Pairings and homotopy (co-)limits
\endheading
\subheading{A.1 Spectra and homotopy (co-)limits}
Recall \cite{A},\cite{T1} that a simplicial {\it prespectrum} $X$ is a
sequence $X_n, n\in {\Bbb N}$ of pointed simplicial sets, together
with maps $ SX_n=S^1\wedge X_n\to X_{n+1}$.  Its homotopy groups are
defined, for all $m\in {\Bbb N}$, as the inductive limit $${\pi}_m(X)
= \text{lim}_n {\pi}_{m+n}(X_n).$$ By $S^1$ we mean the simplicial
circle obtained by identifying the two vertices of the standard one
simplex.  Note that the structure maps are adjoint to maps $X_n\to
\Omega X_{n+1}$.  We say that $X.$ is a {\it spectrum} if these maps
are weak homotopy equivalences for all $n$ and the $X_n$ are all
fibrant. We assume furthermore that all spectra are cofibrant as
prespectra, i.e. that the maps $SX_n=S^1\wedge X_n\to X_{n+1}$ are
inclusions.

 We shall also need smash products of (pre-)spectra. These are only
fully developed in the existing literature for topological spectra;
however we are using simplicial spectra, because there is a fully
developed theory of homotopy limits in the simplicial situation.
Fortunately we shall need only elementary properties of the smash
product, in particular we make no use of associativity and
commutativity of the smash product. We therefore will use the
handicrafted smash product of Boardman, as described in the book of
Adams \cite{A}, replacing the topological spectra in that book by
simplicial prespectra.  On passing to the geometric realization we get
the smash product of Boardman on topological spectra.

Let $\bold I$ be a small category, and $A^*:i\mapsto A^i$, a
contravariant functor from $\bold I$ to the category of spectra.  Then
the homotopy limit of the functor $A^*$ is the spectrum defined,
following \cite{B-K} XI 3.2, and \cite{T1}5.6 as $$
\text{holim}_i(A^i):=
\text{Hom}_{{\bold I}^{\text{op}}} (B({\bold I}\backslash -), A^* ),$$
where $\text{Hom}_{{\bold I}^{\text{op}}}$ means morphisms of
contravariant functors on ${\bold I}$.  Similarily we can define, as
in \cite{B-K} XII 2.1, and \cite{T1} 5.6 and 5.10 the homotopy colimit
$\text{hocolim}_i(A_i )$ of a covariant functor $A.:{\bold I} \to
\text{\bf {Prespectra}}$, to be the difference cokernel of: $$
\coprod_{i\to j} B({\bold I}\backslash j)\triangleright A_i
{\to \atop \to}
\coprod_{k} B({\bold I}\backslash k)\triangleright A_k
$$ in the category of prespectra. Here $X\triangleright Y$ for $X$ a
simplicial set and $Y$ a prespectrum denotes the prespectrum
$X_+\wedge Y$, where $X_+$ denotes the simplicial set $X$ equipped
with a disjoint basepoint.

\bigskip \subheading{A.2 Pairings} \smallskip \subheading{A.2.1}
The following discussion is part of the general theory of homotopy
coends; we make no claim of originality. See for example
\cite{H-V} for a discussion in the case of diagrams of spaces.
Let $A^*$ be an ${\bold I}^{\text{op}}$-diagram of spectra, and $B_*$
and $C_*$ be ${\bold I}$-diagrams of prespectra.  A {\it pairing}
$$A^* \wedge B_* \to C_*$$ consists of pairings $$A^i \wedge B_i \to
C_i $$ for all $i$, such that for all maps $f:i \to j$ in $\bold I$ we
have a commutative diagram: $$ \CD A^j\wedge B_i @>{A(f)\wedge
\text{Id}}>> A^i\wedge B_i @>>> C_i \\ @| @.  @VV{C(f)}V \\ A^j\wedge
B_i @>{\text{Id}\wedge B(f)}>> A^j\wedge B_j @>>> C_j\ .
\endCD $$

\proclaim{Proposition 9} With the notation of the preceeding
definition, a pairing of $\bold I$-diagrams $A^* \wedge B_* \to C_* $
gives rise to a pairing $$ \text{holim}_i(A^i )\wedge
\text{hocolim}_i(B_i ) \to \text{hocolim}_i(C_i )\ . $$ \endproclaim
\demo{Proof} We start by considering a slight generalization of the
homotopy colimit.  Given an ${\bold I}^{\text{op}}$-diagram of pointed
spaces or prespectra $X^* $ and an ${\bold I}$-diagram of pointed
spaces or prespectra $Y_* $, we can form the pointed space or
prespectrum $X^* \dashv Y_* $ defined as the difference cokernel: $$
\coprod_{f:i\to j} X^j\wedge Y_i \underset {\text{Id}\wedge Y(f)}\to
{\overset { X(f)\wedge \text{Id} }\to{{\to \atop \to}}} \coprod_{k}
X^k\wedge Y_k
\quad .
$$ In particular the homotopy colimit of a diagram of prespectra $Y_*
$ may be defined as $$B({\bold I}\backslash *)_+\dashv Y_* \quad .$$
Hence the theorem follows from the more general statement that given a
pairing $\mu : A^*\wedge B_*\to C_*$ as above, there exists for any
$I^{\text{op}}$-diagram $Z^*$ of pointed simplicial sets, a pairing of
prespectra $$ {\text{\bf Hom}}_{I^{\text{op}}}(Z^*,A^*)\wedge
(Z^*\dashv B_*)\to (Z^*\dashv C_*) \quad .  $$

 Given a simplex $\phi\in {\text{\bf Hom}}_{I^{\text{op}}}(Z^*,A^*) $,
for each $i\in I$ we write $\phi_i\in {\text{\bf
Hom}}_{I^{\text{op}}}(Z^i,A^i)$ for its projection into the $i$-th
factor of the product.  Then for each $i$ we have a map $$
\gamma_i :
 {\text{\bf Hom}}_{I^{\text{op}}}(Z^*,A^*)\wedge (Z^i\wedge B_i) \to
(Z^i\wedge C_i) $$ $$ \phi\wedge z^i \wedge b_i \mapsto z^i\wedge
\mu_i (\phi_i(z^i)\wedge b_i), $$ where $\mu_i:A^i\wedge B_i \to C_i$
is the pairing $\mu$ evaluated at $i\in I$. We must verify that these
pairings induce a map between the difference cokernels in the
construction of $Z^*\dashv B*$ and $Z^*\dashv C*$.  That is, given a
map $f:k\to l$ in $\bold I$, we need to know that the following
diagrams commute.  $$ \CD {\text{\bf
Hom}}_{I^{\text{op}}}(Z^*,A^*)\wedge (Z^l\wedge B_k) @>1\wedge
Z(f)\wedge 1 >> {\text{\bf Hom}}_{I^{\text{op}}}(Z^*,A^*)\wedge
(Z^k\wedge B_k) \\ @V\beta_{l,k}VV @VV\gamma_kV \\ Z^l\wedge C_k
@>>Z(f)\wedge 1> Z^k\wedge C_k\ , \endCD $$ and $$
\CD {\text{\bf Hom}}_{I^{\text{op}}}(Z^*,A^*)
\wedge (Z^l\wedge B_k) @>1\wedge
1\wedge B(f) >> {\text{\bf Hom}}_{I^{\text{op}}}(Z^*,A^*)\wedge
(Z^l\wedge B_l) \\ @V\beta_{l,k}VV @VV\gamma_lV \\ Z^l\wedge C_k
@>>1\wedge C(f)> Z^l\wedge C_l\ , \endCD $$ where $$\beta_{k,l} :
\phi\wedge z^l\wedge b_k \mapsto z^l\wedge
\mu_k(A(f)(\phi_l(z^l))\wedge b_k) \quad .  $$ For the first diagram,
we have $$
\align
\gamma_k((1\wedge Z(f)\wedge 1)(\phi\wedge z^l\wedge b_k))
   &= Z(f)(z^l)\wedge \mu_k(\phi_k(Z(f)(z^l)\wedge b_k) \\ &=
Z(f)(z^l)\wedge \mu_k(A(f)(\phi_l(z^l))\wedge b_k) \\ &\quad
\text{(since}\quad A(f)\cdot \phi_l=\phi_k\cdot Z(f)\text{)}\\ &=
(Z(f)\wedge 1)(\beta_{k,l}(\phi\wedge z^l\wedge b_k))\quad .
\endalign $$ While for the second square $$
\align
\gamma((1\wedge 1\wedge B(f))(\phi\wedge z^l\wedge b_k))
   +&= \gamma(\phi\wedge z^l\wedge B(f)(b_k))\\ &= z^l\wedge \mu_l
(\phi_l ( z^l)\wedge B(f)(b_k))\\ &= z^l\wedge C(f) (\mu_k
(A(f)(\phi_l(z^l))\wedge b_k)\\ &\quad (\text{by the ``projection
formula'' for the pairing}\\ &\quad A^*\wedge B_*\to C_* ) \\ &=
(1\wedge C(f))(\beta_{k,l}(\phi\wedge z^l\wedge b_k))\quad , \endalign
$$ and we are done.  \qed \enddemo

\bigskip \subheading{A.2.2}
Suppose we have maps of prespectra $A \overset {f}\to \rightarrow B$
and $E \overset{h} \to \rightarrow F$, and of spectra $C \overset {g}
\to
\leftarrow D$, and pairings
$$ \phi :C \wedge A \rightarrow E $$ and $$ \psi :D \wedge B
\rightarrow F $$ such that the following diagram commutes $$\CD
D\wedge A @>{g\wedge \text{Id}}>> C\wedge A @>>> E \\
@VV{\text{Id}\wedge f}V @.  @VVhV \\ D\wedge B @>>> @>>> F\ .  \\
\endCD\leqno{A.1} $$

We may then define as follows a natural pairing

$$
\mu :{\text {Fiber}} (g) \wedge {\text {Cofiber}} (f) \rightarrow F \ .
$$ Note that ${\text {Fiber}} (g)$ is a spectrum, while ${\text
{Cofiber}} (f) $ is a prespectrum.

\medskip

 For $C$ and $D$ spaces rather than spectra, the homotopy fiber
${\text {Fiber}} (g)$ of $g$ is, by definition, the subcomplex of the
product $D\times (C,*)^{(\Delta[1],*)}$ which makes the square

$$
\CD
 {\text {Fiber}} (g) @>>> (C,*)^{(\Delta[1],*)}\\ @VVV @VVV \\ D @>>>
C
\endCD
$$ become cartesian.  Here $(C,*)^{(\Delta[1],*)}$ is the function
space of pointed maps.  Since $C$ is supposed fibrant, one knows that
the evaluation map
$(C,*)^{(\Delta[1],*)}\to C$ is a Kan fibration. It follows that the
same is true for the map $\text {Fiber} (g)\to D$ induced by
projection onto the first factor. Hence $\text {Fiber} (g)$ is
fibrant.  For $C$ and $D$ spectra, we apply the above construction
degreewise; since $C$ and $D$ are fibrant spectra, the same is true
for $\text {Fiber} (g)$. On the other hand, the mapping cone ${\text
{Cofiber}} (f)$ is the disjoint union of $A
\wedge\Delta[1]$ with $B$, modulo the identification $(a,1) = f(a)$;
note that this is only a prespectrum even if $A$ and $B$ are spectra.
The pairing $\mu$ is defined as follows, we have $$\text
{Cofiber}(f)=B\cup (\Delta[1]\wedge A)\quad .$$ We therefore define
the pairing separately on two pieces: $$\mu :{\text {Fiber}} (g)\wedge
B\to F$$ is induced by the projection ${\text {Fiber}} (g)\to D$ and
the product $\psi:D\wedge B\to F$, while $$\mu :{\text {Fiber}}
(g)\wedge \Delta[1]\wedge A\to F$$ is induced by the adjunction map $$
{\text {Fiber}} (g)\wedge \Delta[1]\to C^{\Delta[1]}\wedge\Delta[1]\to
C $$ followed by the map $h\circ\phi:C\wedge A\to E$.  Using the
commutativity of diagram (A1), one may then check that on the
intersection $CA\cap B\subset {\text {Cofiber}}(f)$
these two maps
agree.

\enddocument